\documentclass[12pt]{article}
\usepackage{graphicx,amsmath,amsthm,amssymb,dsfont, mathrsfs,MnSymbol}
\usepackage{setspace,bbm}
\def\ZZ{\mathbb Z}
\def\PP{\mathbb P}

\def\NN{\mathbb N}
\def\QQ{\mathbb Q}

\def\DD{\mathbb D}
\def\FF{\mathbb F}

\def\VV{\mathbb V}

\def\SS{\mathbb S}
\def\EE{\mathbb E}

\def\cE{\mathcal E}

\def\sL{\mathscr L}
\def\cL{\mathcal L}

\def\cM{\mathcal M}

\def\cJ{\mathcal J}

\def\cF{\mathcal F}
\def\ff{\mathfrak f}
\def\fg{\mathfrak g}

\def\fa{\mathfrak a}
\def\fb{\mathfrak b}
\def\fc{\mathfrak c}
\def\fm{\mathfrak m}

\def\cE{\mathcal E}

\def\fd{\mathfrak d}

\def\cP{\mathcal P}

\def\cN{\mathcal N}

\def\fT{\mathfrak T}

\def\bA{{\bf A}}

\def\bT{{\bf T}}
\def\mt{{\middlebar{t}}}
\def\bx{{\bf x}}

\def\bbx{{\mathbbm x}}

\def\by{{\bf y}}

\def\kk{{\mathbbm k}}

\def\bb{{\bf b}}

\def\brho{\boldsymbol{\rho}}

\def\bo{{\bf 0}}
\def\b1{{\bf 1}}
\def\d1{\mathds{ 1}}
\def\mod{{\; \rm mod \;}}
\def\rd{{\rm d}}
\def\ad{{\rm and}}

\def\with{{\rm with}}
\def\where{{\rm where}}
\def\for{{\rm for}}

\def\vol{{\rm vol}}

\def\qed{ \ \vrule width.2cm height.2cm depth0cm\smallskip}

\begin{document}


\title{Central Limit Theorem for $(t,s)$-sequences, I}
\author{Mordechay B. Levin}

\date{}

\maketitle

\begin{abstract}
Let $ (X_n)_{n \geq 0} $ be a digital $(t,s)$-sequence in base $2$, $\cP_{m} =(X_n)_{n=0}^{2^m-1} $, and let $D( \cP_{m}, Y  )$
 be the local discrepancy of  $ \cP_{m}$. Let $T \oplus Y$ be the digital addition of $T$ and $Y$, and let
 $$
           \cM_{s,p} (\cP_{m})  =\Big( \int_{[0,1)^{2s}} |D( \cP_{m} \oplus T , Y  ) |^p \rd T  \rd Y \Big)^{1/p} .
 $$
In this paper, we prove that  $D( \cP_{m} \oplus T , Y  ) / \cM_{s,2} (\cP_{m})$ weakly converge to the standard Gaussisian distribution for $m \rightarrow \infty$, where $T,Y$ are uniformly distributed random variables in $[0,1)^s$. In addition, we prove that
\begin{equation} \nonumber
     \cM_{s,p} (\cP_{m}) / \cM_{s,2} (\cP_{m})   \to  \frac{1}{\sqrt{2\pi}}\int_{-\infty}^{\infty} |u|^p e^{-u^2/2} \rd u        \quad {\rm for} \; \; m \to  \infty , \;\; p>0.
\end{equation}

\end{abstract}
%
{\bf 1. \;\; Introduction }\\
%
 Let $(\beta_{n})_{n = 0}^{N-1}$ be an $N$-element point set in the
 $s$-dimensional unit cube $[0,1)^s$. The local {\bf discrepancy} function of $(\beta_{n})_{n = 0}^{N-1}$ is defined as
\begin{equation} \nonumber 
   D( (\beta_{n})_{n = 0}^{N-1}, Y  )= \sum\nolimits_{n=0}^{N-1}   \b1_{B_{Y}}(\beta_{n}) - N y_1 \cdots y_s,
\end{equation}
where $Y=(y_1,...,y_s)$, $B_{Y}=[0,y_1) \times \cdots \times [0,y_s) $,
 $\b1_{B_{Y}}(X) =1, \; {\rm if} \; X  \in B_{Y}$,
and $   \b1_{B_{Y}}(X) =0,$  if $
X   \notin B_{Y}$.
We define the $\emph{L}_p$  discrepancy of
$(\beta_{n})_{n = 0}^{N-1}$ as
\begin{equation}   \nonumber 
       D_{\infty} ((\beta_{n})_{n = 0}^{N-1}) =
    \sup_{ 0<y_1, \ldots , y_s \leq 1} \; |
    D((\beta_{n})_{n = 0}^{N-1},Y) |,\;\;  D_{p}((\beta_{n})_{n = 0}^{N-1})=\left\|   D((\beta_{n})_{n = 0}^{N-1},Y) \right\|_{p}  ,
\end{equation}
\begin{equation*} 
     \left\| f(Y) \right\|_{q}= \Big( \int_{[0,1]^s}  |f(Y)|^p d Y \Big)^{1/p}.
\end{equation*}
{\bf Definition 1.}
 A sequence $(\beta_n)_{n\geq 0}$ is of {\it
low discrepancy} (l.d.s.) if  $D_{\infty}
((\beta_n)_{n=0}^{N-1})=O(\log^s N) $.
A point set $(\beta_{n,N})_{n=0}^{N-1}$ is of
 low discrepancy (l.d.p.s.) if $ D_{\infty} ((\beta_{n,N})_{n=0}^{N-1})=O(\log^{s-1}
N)$.
 For examples of such a sequence, see, e.g., [BC], [DiPi1], and  [Ni].\\

A subinterval $U$ of $[0,1)^s$  of the form
$$  U = \prod_{i=1}^s [a_ib^{-d_i},(a_i+1)b^{-d_i}),   $$
 with $a_i,d_i \in \ZZ, \; d_i \ge 0, \; 0  \le a_i < b^{d_i}$ for $1 \le i \le s$ is called an
{\it elementary interval in base $b \geq 2$}.\\ 

 {\bf Definition 2}. {\it Let $0 \le t \le m$  be  integers. A {\sf $(\mt,m,s)$-{\sf net in base $b$}} is a point set
$\bx_0,...,\bx_{b^m-1}$ in $ [0,1)^s $  such that $\# \{ n \in [0,b^m -1] | \bx_n \in U \}=b^t$   for every elementary interval $U$ in base  $b$ with
$\vol(U)=b^{t-m}$.}\\  \\

%

First constructions of dyadic  $(\mt,m,s)$ nets and $(\mt,s)$ sequences  were given by Sobol [So]. For other
 constructions and references see [DiPi1] and [Ni]. These nets and  sequences are of low discrepancy (see, e.g., [Ni, p. 56,60]).

It is known that
\begin{equation} \nonumber 
\emph{D}_{s,p}((\beta_{n,N})_{n=0}^{N-1})>C_{s,p}(\log N)^{\frac{s-1}{2}} \;\;\;
\end{equation}
for all $N$-point sets $(\beta_{n,N})_{n=0}^{N-1}$  with some $C_{s,p} >0$ (see  Roth for $p=2$, Schmidt  for $p>1$ [BC].\\

{\bf Definition 3.}  A sequence  of point sets $((\beta_{n,N})_{n=0}^{N-1})_{N=1}^{\infty}$ is of
 $L_p$ low discrepancy (l.d.p.s.) if $ D_{s,p} ((\beta_{n,N})_{n=0}^{N-1})=O((\log
N)^{(s-1)/2}) $ for $ N \rightarrow \infty $.

The existence of $L_p$ l.d.p.s. was proved by Roth for $p=2$ and by Chen  for $p>1$ [Ch].
The first explicit construction of $L_p$ l.d.p.s. was obtained by Chen and Skriganov for $p=2$ and by Skriganov  for $p>1$ (see [ChSk], [Skr2]).
The next explicit construction of $L_p$ l.d.p.s. was proposed by Dick and Pillichshammer
 (see [Di], [DiPi2],  [Ma]).

The first lower discrepancy bound for $p \in (0,1]$ was obtained by Skriganov [Skr1]: \\

We write $\NN$ for the set of all positive integers, $\NN_0$ for the set of all non-negative integers.
For $i \in \NN_0$, we put
\begin{equation} \nonumber
\QQ(2^i) =\Big\{\frac{j}{2^{i}}\;|\;j=0,1,..., 2^i-1 \Big\} \;\;\ad \;\;\QQ^s(2^i)= \{(x_1,...,x_s) | x_j \in \QQ(2^i),  j \in [1,s]\}.
\end{equation}
The points of $\cup_{i \geq 0} \QQ^s(2^i)$ are called {\it dyadic rational points}.
Any $y \in [0,1)$ can be represented in the form
\begin{equation} \label{In5}
 y= \sum_{a \geq 1} y_a 2^{-a}  , \qquad \where \quad y_a\in \{ 0,1\}, \; a \in \NN.
\end{equation}
For any two points $x$ and $y$ in $[0,1)$, we define their sum $x \oplus y$ by
\begin{equation} \label{In6}
 (x \oplus y)_a \equiv x_a + y_a\mod 2, \qquad  a \in \NN.
\end{equation}
For dyadic rational numbers, we always use the finite expansion. In this paper, we will use $\oplus$ only for the case that $y$ (or $x$) is a dyadic rational number. For this case \eqref{In6} define the addition $\oplus$ in a proper sense. For vectors $\bx, \by \in [0,1)^s$, we use the notation $\bx \oplus \by$ to denote the component-wise addition $\oplus$.

Let
\begin{equation} \nonumber 
 \cM_{m,p}(\cP_{m}) = \Big(  2^{-sm} \sum_{T \in \QQ^s(2^m)  }  D_p(\cP_{m} \oplus T)^p  \Big)^{1/p},   \qquad 0 <p< \infty.
\end{equation}
For dyadic $(\mt,m,s)$-net $\cP_m$,
 Skriganov [Skr3, Theorem 2.1,\; Theorem 2.2 ] proved  for $p>0$ that
\begin{equation} \label{In8}
   m^{-(s-1)/2}  \cM_{s,p} (\cP_{m}) \in [2^{-2s-2}(1+1/p) s^{-(s-1)/2}, sp2^{\mt+2}].
\end{equation}
The case of $L_2$-discrepancy was studied earlier (see, e.g., [DiPi1, Sections 16.5 and 16.6]).

In this paper, we will make \eqref{In8} precise  for the case of digital $(\mt,m,s)$-sequences in base $2$.
We will use the following definitions of  digital $(\mt,m,s)$ nets, digital  $(\mt,s)$ sequences  and
$(\bT,s)$ sequences: \\ 

 {\bf Definition 4.} (\cite[\S 4.4]{DiPi1}) { \it
Let $m, s \geq 1$ be integers. Let $C^{(1,m)},...,$ $C^{(s,m)}$ be $m \times m$ matrices over $\FF_2$.
Now we construct $2^m$ points in $[0, 1)^s$.
 For $ n= 0, 1,...,2^m-1$, let $n =\sum^{m-1}_{j=0} e_j(n) 2^{j}$
be the dyadic expansion of $n$.  For $r = 0,1,...$,
we choose bijections  $\upsilon_r : \ZZ_b \to  \FF_b $ with  $\upsilon_r(0) = 0$, and for $i = 1, 2, . . . , s$
and $j = 1, 2, . . . $ we choose bijections $\varsigma_{i,j} : \FF_2 \to  \ZZ_b $.
We map the vectors
\begin{equation} \nonumber 
	y_{n,i}=(y_{n,i,1},...,y_{n,i,m}),\quad y_{n,i,j}=\sum_{r=0}^{m-1} \upsilon_r(e_r(n)) c^{(i,m)}_{j,r}\in  \FF_2
\end{equation}
to the real numbers
\begin{equation}  \nonumber 
   x_{n,i} =\sum_{j=1}^m \varsigma_{i,j} (y_{n,i,j})/2^j
\end{equation}
to obtain the point
\begin{equation} \nonumber
   \bx_n= (x_{n,1},...,x_{n,1}) \in [0,1)^s.
\end{equation}

The point set  $ \{\bx_0,...,\bx_{2^m-1} \}$ is called a dyadic {\sf digital net}  (with {\sf generating matrices} $(C^{(1,m)},...,C^{(s,m)}) $).

For $m = \infty$, we obtain a sequence $\bx_0, \bx_1,...$ of   points in $[0, 1)^s$  which is called a dyadic  {\sf digital sequence} $($with {\sf generating matrices} $(C^{(1,\infty)},...,C^{(s,\infty)}) )$.}

We abbreviate  $C^{(i,m)}$ as $C^{(i)}$ for $m \in \NN$ and for $m=\infty$. \\ 

 {\bf Definition 5.} (\cite[Definition 4.30]{DiPi1}) { \it For a given dimension $s \geq 1$,  and a
function $\bT : \NN_0 \to \NN_0$ with $\bT(m) \leq  m$ for all $m \in \NN_0$, a sequence $(\bx_0,\bx_1, . . .)$
of points in $[0, 1)^s$ is called a dyadic $(\bT, s)$-sequence  if for all integers $m \geq 1$
and $k \geq  0$, the point set consisting of the points $x_{k2^m}, . . . ,x_{k2^m+2^m-1}$ forms
a dyadic $(\bT(m),m, s)$-net. }\\

In this paper, we will prove \\ \\
{\bf Theorem 1.} {\it Let  $(\cP_m)_{m \geq 1}$ be a sequence of  dyadic digital $(\mt_m,m,s)$-nets with
 $\mt_m\leq  1/10 \; \log_2 m$, $s\geq 2$. Then}
\begin{equation}   \nonumber 
   2^{-sm} \sum_{T \in \QQ^s(2^m)} \vol  \Big\{ Y \in [0,1)^s \; : \;\frac{ D( \cP_m \oplus T,Y)}{ \cM_{s,2} (\cP_{m})} <w    \Big\}
	\stackrel{m \to \infty }{\longrightarrow} \Phi(w) ,
\end{equation}
where $\Phi(w) = \frac{1}{\sqrt{2\pi}}\int_{-\infty}^{w} |u|^p e^{-u^2/2} \rd u  $. \\ \\
{\bf Theorem 2.} {\it Let  $(\cP_m)_{m \geq 1}$ be a sequence of  dyadic digital $(\mt,m,s)$-nets
 and  $s\geq 2$. Then}
\begin{equation}  \nonumber 
  \frac{ \cM_{s,p} (\cP_{m}) ) }{ \cM_{s,2} (\cP_{m})}  \stackrel{m \to \infty }{\longrightarrow}  \chi_p, \quad \where \quad
  \chi_p =   \frac{1}{\sqrt{2\pi}}\int_{-\infty}^{\infty} |u|^p e^{-u^2/2} \rd u, \;\;\; p>0
\end{equation}
and  $\chi_{2r} =(2r)!/(2^{r}r!) $ for  integers  $r \geq 1$. \\ \\
{\bf Corollary 1.} {\it Let  $(X_n)_{n\geq 0}$ be a  dyadic digital $(\bT,s)$-sequence  with \\
 $\bT(m) \leq  1/10 \; \log_2 m$, $\dot{\cP}_m =(X_n)_{n= 0}^{2^m-1} $, $s\geq 2$ and $\ddot{\cP}_m =(n/2^{m}, X_n)_{n= 0}^{2^m-1} $. Then}
\begin{equation} \nonumber 
   2^{-sm} \sum_{T \in \QQ^s(2^m)} \vol \Big\{  Y \in [0,1)^s \; : \; \frac{ D( \dot{\cP}_m \oplus T,Y)}{ \cM_{s,p} (\dot{\cP}_m)} <w    \Big\}
	\stackrel{m \to \infty }{\longrightarrow} \Phi(w) ,
\end{equation}
\begin{equation}  \nonumber 
   2^{-(s+1)m} \sum_{T \in \QQ^{s+1}(2^m)} \vol \Big\{  Y \in [0,1)^{s+1} \; : \; \frac{ D( \ddot{\cP}_m \oplus T,Y)  }{ \cM_{s+1,p} (\ddot{\cP}_m)} <w    \Big\}
	\stackrel{m \to \infty }{\longrightarrow} \Phi(w) ,
\end{equation}
and
\begin{equation} \nonumber
  \frac{ \cM_{s,p} (\dot{\cP}_m) ) }{ \cM_{s,2} (\dot{\cP}_m)}  \stackrel{m \to \infty }{\longrightarrow}  \chi_p,
  \qquad   \frac{ \cM_{s+1,p} (\ddot{\cP}_m)  }{ \cM_{s+1,2} (\ddot{\cP}_m)}  \stackrel{m \to \infty }{\longrightarrow}  \chi_p, \;\;\; \for \; \limsup_{m \to \infty} \bT(m) < \infty, \; p>0.
\end{equation} \\

In a forthcoming paper, we will prove Theorem 1 and Theorem 2 for the case of the arbitrary base.
In another forthcoming paper, we will consider 2-order digital $(t,s)$-sequences (see the definition in [Di]) and we will prove the Central Limit Theorem (CLT) in Theorem 1 for an  $s$-dimensional random variable $Y$ instead of an $2s$-dimensional random variable $(T,Y)$.\\

Now we describe the structure of the paper.
\S2 and \S3 are auxiliary chapters.
In \S4, we get four moments estimates of the discrepancy function. We calculate the four moments in the proper sense and the Levi conditional expectation.
In \S5, we use a variant of martingale CLT to prove Theorem 1. Theorem 2 is a simple corollary
of Theorem 1 and \eqref{In8}.
\\ \\
{\bf 2. \;\; Skriganov's formula for discrepancies.}

We will use notations from [Skr3].
In the one-dimensional case,  the Rademacher functions $r_a(y)$, $y \in [0,1)$, $a \in \NN$, can be defined by
\begin{equation} \nonumber
   r_a(y) =  1 - 2 y_a,
\end{equation}
where $y_a$ are the coefficients in the dyadic expansion \eqref{In5}. It is convenient to put
$r_0(y) \equiv 1 $.
The $s$-dimensional Rademacher functions $r_A(Y)$, $Y=(y_1,...,y_s) \in [0,1)^s$, $A=(a_1,...,a_s) \in \NN^s_0$, are defined by
\begin{equation} \nonumber
   r_A(Y) = \prod_{i=1}^s r_{a_i}(y_i).
\end{equation}

Let $y_i =0.y_{i,1}y_{i,2}...=\sum_{j \geq 1}  y_{i,j} 2^{-j}$, with $y_{i,j} \in \{0,1\}$,
 $i=1,...,s$.
We define the truncation
\begin{equation}  \label{Le2-1}
        y_i^{(m)} =\sum_{1 \leq j \leq m} y_{i,j} 2^{-j} .
\end{equation}
If $Y = (y_1, . . . , y_s)  \in [0, 1)^s$, then the truncation $ Y^{(m)}$ is defined coordinatewise, that is $Y^{(m)}=
( y_1^{(m)}, . . . , y_s^{(m)})$. Let
\begin{equation} \nonumber
   \delta_2(a) =   \begin{cases}
    1,  & \; {\rm if}  \;  a \equiv 0  \mod 2,\\
    0, &{\rm otherwise}
  \end{cases}  \qquad \ad \qquad
   \d1(\fT) =   \begin{cases}
    1,  & \; {\rm if}  \;  \fT  \;{\rm is \;true},\\
    0, &{\rm otherwise}
  \end{cases}.
\end{equation}
It is known that
\begin{equation} \label{Le2-3}
  \delta_2(a) = \sum_{\ell =0,1} e(\ell a), \quad \where \quad e(a) =exp(\pi i a).
\end{equation}
Consider the elementary intervals
\begin{equation} \nonumber
   \Pi_a =  [ 2^{-a}, 2^{1-a}), \qquad a \in \NN.
\end{equation}
It is convinient to put $\Pi_0 = [0,1)$.
Introduce elementary boxes of the form
\begin{equation} \label{Le2-4}
   \Pi_A =\Pi_{a_1} \times \cdots \times \Pi_{a_s}, \qquad A=(a_1,...,a_s) \in \NN_0^s.
\end{equation}
Each such box has volume $\vol (\Pi_A ) = 2^{-a_1-\cdots -a_s}$.
Let
\begin{equation} \nonumber
   \lambda_A(Y) =\b1_{\Pi_{A}} (Y) -\vol( \Pi_A) .
\end{equation}
We put
\begin{multline} \label{Le2-5}
E^{(m)}_{T} ( f(T,Y)) = \frac{1}{2^{sm}} \sum_{\substack{ t_{i,j} \in \{0,1\} \\ 1 \leq i \leq s,\; 1 \leq j \leq m}} f(T^{(m)}, Y^{(m)}), \\
E^{(m)}_{Y} ( f(T,Y)) = \frac{1}{2^{sm}} \sum_{\substack{ y_{i,j}\in \{0,1\} \\ 1 \leq i \leq s,\; 1 \leq j \leq m}} f(T^{(m)}, Y^{(m)}),  \\
           E^{(m)}_{Y,T} ( f(T,Y)) = E^{(m)}_{T} (  E^{(m)}_{Y} ( f(T,Y) ) ), \\
            \EE_{Y,T,m} ( f(T,Y)) =  \frac{1}{2^{sm}} \sum_{\substack{ t_{i,j}\in \{0,1\} \\ 1 \leq i \leq s,\; 1 \leq j \leq m}} \int_{[0,1)^s} f(T, Y) \rd Y.
\end{multline}
It is easy to see that
\begin{equation} \label{Le2-5a}
   E^{(m)}_{Y,T}  ( f(T,Y)) =  \EE_{Y,T,m} ( f(T^{(m)},Y^{(m)})). 
\end{equation}
In the following, we mean by $\PP_{m}=(\bbx_n)_{n=0}^{2^m-1}$ a dyadic digital $(\mt,m,s)$-net. \\ \\
We define the {\it micro-local discrepancy} by
\begin{multline}\label{Le2-7}
   \lambda_A ( \PP_{m}^{(m)} \oplus Y^{(m)}  ) = \sum_{\bbx \in \PP_{m}}
	\lambda_A ( \bbx^{(m)} \oplus Y^{(m)}  ) \\
   = 	\sum_{\bbx \in \PP_{m}}
	\big( \b1_{\Pi_A} ( \bbx^{(m)} \oplus Y^{(m)}  ) - \vol (\Pi_A) \big).
\end{multline}
Using  [Skr3, Lemma 4.3, Lemma 6.1, and ref. (4.6), (4.30), (5.8)], we obtain : \\ \\
{\bf Lemma A.} {\it For each $m \in \NN$, the local discrepancy $ D( \PP_{m}, Y  )$ has the representation }
\begin{equation}  \nonumber 
   D( \PP_{m}, Y  )=  D^{(m)}( \PP_{m}, Y  ) + \cE^{(m)}( \PP_{m}, Y  ),
\end{equation}
with
\begin{equation}  \nonumber 
   D^{(m)}( \PP_{m}, Y  )  = 2^{-s} \sum_{A \in I_m^s} (-1)^{\kappa(A)} \lambda_A ( \PP_{m}^{(m)} \oplus Y^{(m)}  )
   r_A(Y)
\end{equation}
and
 \begin{equation}    \nonumber 
         \EE_{Y,T,m} \big( |\cE^{(m)}(\PP_{m}\oplus T , Y ) |^2 \big)      \leq (s 2^{\mt})^2,
\end{equation}
where $\kappa(A)$ is  the number of non-zero elements in $A$ and $I_m =\{0,1,...,m\}$.\\  \\ \\
{\bf 3. Auxiliary lemmas.}

 We will use notation $A \ll B$  equal to $A = O(B)$.
Let
\begin{equation} \nonumber
  K_m  =\Big\{  \fm = (\fm_1,...,\fm_s)\;
  :   \fm_i=(\fm_{i,1},...,\fm_{i,m}), \;
   \;\fm_{i,j} \in \FF_2, \;\;
  j \in [1,  m], i \in [1,s]\Big\},
\end{equation}	
 $K_m^{*} =  K_m \setminus \{ \bo\}$, and let
\begin{equation}  \nonumber 
     \PP^{\bot}_m =\Big\{  \fm \in K_m \; : \;  \sum_{i=1}^s \sum_{j=1}^m \fm_{i,j} \bbx_{n,i,j}  =0 \;\; \forall \;\;
						n \in [0,2^m-1] \Big\}, \quad  \;\;  \PP^{\bot,*}_m   =  \PP^{\bot}_m  \setminus \{ \bo\}.
\end{equation}
For  any vector $\bb = (b_1,...,b_m) \in \FF^m_2$, let
\begin{equation} \nonumber 
 \rho(\bb) = 0 \;\; {\rm if} \;\; \bb = \bo \quad
 \ad  \quad  \rho(\bb) = \max \{j \; : \; b_j \neq 0  \}\;\; {\rm if} \;\; \bb \neq \bo,
\end{equation}
and let
\begin{equation}  \nonumber 
     \brho(\fb) = \sum_{1 \leq i \leq s}  \rho(\bb_i), \quad \for \quad \fb = (\bb_1,...,\bb_s) \in \FF^{ms}_2.
\end{equation}
We put
\begin{equation} \nonumber
        \brho(\PP^{\bot}_m)  =\min \{ \brho(\fm) \;  : \; \fm \in \PP^{\bot}_{m} \setminus \{ \bo\}  \}.
\end{equation} \\
{\bf Lemma B}. (\cite[Theorem 7.5]{DiPi1},  [Skr1, Theorem 4.2])  {\it The net $\PP_{m}=(\bbx_n)_{n=0}^{2^m-1}$ is  $(\mt,m,s)$ dyadic digital if and only if $ \brho(\PP_m^{\bot}) \geq  m -\mt +1$.} \\ \\
{\bf Lemma C}. (  [Skr2, Lemma 2.2])  {\it Let $\PP_{m}=(\bbx_n)_{n=0}^{2^m-1}$ be a
dyadic digital $(\mt,m,s)$-net, $A=(a_1,...,a_s)$, and let $a_0(A):=a_1+ \cdots +a_s \geq \brho(\cP_m^{\bot}) $. Then} \\
\begin{equation}  \nonumber 
\# G_A \leq 2^{a_0(A) - \brho(\PP_m^{\bot}) +1 }, \quad \with \quad G_A=\{ \fm \in \PP_m^{\bot} \; : \; \rho(\fm_i) \leq a_i, \; i\in [1,s]  \}.
\end{equation} \\
{\bf Lemma D}. (  Skr2, Lemma 4.1)  {\it Let $\PP_{m}=(\bbx_n)_{n=0}^{2^m-1}$ be a
dyadic digital $(\mt,m,s)$-net, and let $\fm \in K_m$. Then} \\
\begin{equation} \nonumber
      \sum_{n=0}^{2^m-1} e\Big(\sum_{i=1}^s \sum_{j=1}^{m} \fm_{i,j}\bbx_{i,j} \Big)  =
       \begin{cases}
    2^m,  & \; {\rm if}  \;  \fm \in \PP_m^{\bot},\\
    0, &{\rm otherwise}
  \end{cases} .
\end{equation} \\

For each $Z \in \QQ^s(2^m)$, the shift $\PP_{m} \oplus Z$ is also a $(\mt,m,s)$-net, and it follows from \eqref{Le2-7} (see also [Skr3, p.205]) that
 \begin{equation}   \nonumber 
 | \lambda_A ( \PP_{m}^{(m)} \oplus Z )  | \leq 2^{\mt} \quad \ad \quad  \lambda_A ( \PP_{m}^{(m)} \oplus Z )=0 \quad  \quad {\rm if \;\; vol}( \Pi_A) \geq 2^{\mt-m}.
\end{equation}

Let $V_0 = 10 \; \log_2 m$,   $I_m =\{0,1,...,m\}$,
\begin{multline}\label{Le3-5}
  I_{m,s,k} =\big\{A =(a_1,...,a_s)\in I_m^s \; : \; \max_i a_i=k, \; a_0(A) \in ( m-\mt_m,m+V_0) \big\}, \\
I_{m,s,m+1} =\big\{A \in I_m^s \; : \; a_0(A)=a_1+ \cdots +a_s  \geq  m+V_0 \big\}.
\end{multline}
It easy to verify that
\begin{equation}    \label{Le3-6}
  \# I_{m,s,k} \leq s k^{s-2}(t+V_0), \qquad k \in [1,m].
\end{equation} \\ \\
{\bf Lemma 1.}  {\it Let $\cP_{m}=(\bx_n)_{n=0}^{2^m-1}$ be the  dyadic digital   $(\mt_m,m,s)$-net. Then  }
\begin{multline}\label{Le3-7}
     D( \cP_{m} \oplus T, Y  )=  D^{(m)}( \cP_{m} \oplus T, Y  ) + \cE^{(m)}( \cP_{m}\oplus T, Y  ),\\
     D^{(m)}( \cP_{m} \oplus T, Y  )  = \sum_{k=1}^m \DD_k(T,Y) +R(T,Y), \\
\with \qquad \qquad \qquad \qquad \qquad \qquad \qquad \qquad \qquad \qquad \qquad \qquad \qquad\qquad \qquad \qquad \qquad \qquad \qquad \qquad\\
   \DD_k =  \DD_k(T,Y) = \sum_{A \in I_{m,s,k}}  \Psi_A,   \quad R(T,Y) = \sum_{A \in I_{m,s,m+1}}  \Psi_A  \\
     \Psi_A =  (-1)^{\kappa(A)}  2^{m-s -a_1-\cdots -a_s}   r_A(Y) \Lambda_A(T \oplus Y^{(m)}  ),\\
\where \qquad \qquad \qquad \qquad \qquad \qquad \qquad \qquad \qquad \qquad \qquad \qquad \qquad\qquad \qquad \qquad \qquad \qquad \qquad \qquad \\
    \Lambda_A(T \oplus Y^{(m)} ) =
 \sum_{\substack{  \fm \in \cP^{\bot,*}_m, \; \rho(\fm_i) \leq a_i \\ i=1,...,s   }} e\Big(\sum_{i=1}^s \big( \sum_{j=1}^{m} \fm_{i,j}(t_{i,j}+y_{i,j} )+ \fm_{i,a_i} \big) \Big).
\end{multline} \\ \\
{\bf Proof.}
In view of \eqref{In6}, \eqref{Le2-1} and \eqref{Le2-4}, we get
\begin{equation}  \nonumber 
   \b1_{\Pi_A} ( \bx_n \oplus Y^{(m)}  ) = \prod_{i=1}^s \prod_{j=1}^{a_i -1}
	\d1( x_{n,i,j} \oplus y_{i,j} =0) \d1( x_{n,i,a_i} \oplus y_{i,a_i} =1) .
\end{equation}
By \eqref{Le2-3}, we obtain
\begin{equation}    \nonumber
 \d1( x_{n,i,j} \oplus y_{i,j} =0) =2^{-1} \sum_{\fm_{i,j} \in \{0,1\}} e(\fm_{i,j}(x_{n,i,j} + y_{i,j} ))
\end{equation}
and
\begin{equation}    \nonumber
 \b1_{\Pi_A } ( \bx_n  \oplus Y^{(m)}  )  =\frac{1}{2^{a_1+\cdots a_s}}
 \sum_{\substack{  \fm_{i,j} \in \{0,1\} \\ 1 \leq i \leq s,\; 1 \leq j \leq a_i}} e\Big(\sum_{i=1}^s \big( \sum_{j=1}^{a_i} \fm_{i,j}(x_{n,i,j}^{(m)}+ y_{i,j}^{(m)}) + \fm_{i,a_i} \big) \Big).
\end{equation}
From \eqref{Le2-7}, we have
\begin{equation}  \nonumber
   \lambda_A ( \cP_{m} \oplus Y^{(m)}  ) = \frac{1}{2^{a_1+\cdots + a_s}}
 \sum_{\substack{  \fm_{i,j} \in \{0,1\} \\ 1 \leq i \leq s,\; 1 \leq j \leq a_i}} \Xi_{\fm} e\Big(\sum_{i=1}^s \big(\sum_{j=1}^{a_i} \fm_{i,j} y_{i,j}+ \fm_{i,a_i} \big) \Big)  - \frac{2^m}{2^{a_1+\cdots + a_s}}
\end{equation}
with
\begin{equation}  \nonumber
\Xi_{\fm} =
    	\sum_{n=0}^{2^m-1} e\big(\sum_{i=1}^s \sum_{j=1}^{a_i} \fm_{i,j} x_{n,i,j} \big).
\end{equation}
Applying Lemma D, we get
\begin{equation}  \nonumber
\Xi_{\fm} = 2^m \;\; \for \;\;  \fm \in \cP^{\bot}_m,  \quad  \quad \Xi_{\fm} = 0 \;\; \for
\;\;      \fm \notin \cP^{\bot}_m
\end{equation}
and
\begin{multline}  \nonumber 
   \lambda_A ( \cP_{m} \oplus Y^{(m)}  ) = 2^{m -a_1-\cdots -a_s}
 \sum_{ \fm \in \cP^{\bot,*}_m} e\Big(\sum_{i=1}^s \big( \sum_{j=1}^{a_i} \fm_{i,j} y_{i,j}+ \fm_{i,a_i} \big) \Big) \\
 =    2^{m- a_1-\cdots -a_s}  \sum_{\substack{  \fm \in \cP^{\bot,*}_m, \; \rho(\fm_i) \leq a_i \\ i=1,...,s   }} e\Big(\sum_{i=1}^s \big( \sum_{j=1}^{m} \fm_{i,j} y_{i,j}+ \fm_{i,a_i} \big) \Big).
\end{multline}
Using Lemma A and \eqref{Le3-7}, we obtain
\begin{equation}  \nonumber 
   D( \cP_{m} \oplus T, Y  )=  D^{(m)}( \cP_{m} \oplus T, Y  ) + \cE^{(m)}( \cP_{m} \oplus T, Y  ),
\end{equation}
with
\begin{multline} \nonumber 
   D^{(m)}( \cP_{m} \oplus T, Y  )  = 2^{-s} \sum_{A \in I_m^s} (-1)^{\kappa(A)} \lambda_A ( \cP_{m}^{(m)} \oplus T \oplus Y^{(m)}  )
   r_A(Y)  \\
  =  \sum_{A \in I_m^s}  (-1)^{\kappa(A)}
		2^{m-s -a_1-\cdots -a_s}  \Lambda_A(T \oplus Y^{(m)}  )  r_A(Y).
\end{multline}
Bearing in mind that $ \brho(\cP_m^{\bot}) \geq  m -\mt_m +1$, we obtain that $\Lambda_A(T \oplus Y^{(m)}  )=0 $ for $a_0(A) =a_1 + \cdots + a_s \leq m -\mt_m$. By \eqref{Le3-7}, we have
\begin{multline}  \nonumber
    D^{(m)}( \cP_{m} \oplus T, Y  )  =  \sum_{A \in  I_m^s,\; a_0(A) > m -\mt_m }
			 \Psi_A
= \sum_{k=1}^m \DD_m', \;\; \with \;\;  \DD_m'=\sum_{A \in \tilde{I}_{m,s,k} }
 \Psi_A, \\
     \Psi_A =  (-1)^{\kappa(A)}  2^{m-s -a_1-\cdots -a_s}   r_A(Y) \Lambda_A(T \oplus Y^{(m)}  ),
\end{multline}
where $\tilde{I}_{m,s,k} =\big\{A \in I_m^s \; : \; \max_i a_i=k, \; a_0(A) > m-\mt_m \big\}$.\\
Taking into account \eqref{Le3-5}, we get
\begin{equation}  \nonumber
\bigcup_{k=1}^m  \tilde{I}_{m,s,k} =   \bigcup_{k=1}^{m+1}  I_{m,s,k}, \quad \with \quad     I_{m,s,k_1}\cap I_{m,s,k_2} = \emptyset \;\; \for \;\; k_1 \neq k_2    ,
\end{equation}
where
\begin{multline}    \nonumber
  I_{m,s,k} =\big\{A =(a_1,...,a_s)\in I_m^s \; : \; \max_i a_i=k, \; a_0(A) \in ( m-\mt_m,m+V_0) \big\}, \\
I_{m,s,m+1} =\big\{A \in I_m^s \; : \; a_0(A)=a_1+ \cdots +a_s  \geq  m+V_0 \big\}.
\end{multline}
Hence
\begin{multline}    \nonumber
     D^{(m)}( \cP_{m} \oplus T, Y  )  = \sum_{k=1}^m \DD_k(T,Y) +R(T,Y), \\
\with \qquad
     \DD_k(T,Y) = \sum_{A \in I_{m,s,k}}  \Psi_A,  \quad  R(T,Y) = \sum_{A \in I_{m,s,m+1}}  \Psi_A.
\end{multline}
Therefore Lemma 1 is proved. \qed  \\ \\
Let
\begin{equation}  \nonumber %
\fm^{[\ell]} =(\fm^{[\ell]}_1,...,\fm^{[\ell]}_s),  \quad \with \quad  \fm^{[\ell]}_i =(\fm^{[\ell]}_{i,1},...,\fm^{[\ell]}_{i,m}), \quad \ell =1,...,4.
\end{equation} \\ \\
{\bf Lemma 2.}  {\it Let   $\mt_m \leq 1/10 \;  \log_2 m$. Then}
\begin{equation}   \nonumber 
        \EE_{Y,T,m}((\cE^{(m)}(\cP_{m}\oplus T , Y ) +R(T,Y))^2) \ll  m^{1/5}.
\end{equation} \\
{\bf Proof.} By  \eqref{Le3-7}, we obtain
\begin{multline} \nonumber
  \tau:= E_T^{(m)}  R^2(T,Y) = \sum_{A_1, A_2 \in I_{m,s,m+1}} (-1)^{\kappa(A_1)+\kappa(A_2)}
	2^{2m-2s -a_{0}(A_1)- a_{0}(A_2)}    \\
 \times  r_{A_1}(Y)  r_{A_2}(Y) \; \varpi, \quad \where \quad \varpi:=   E_T^{(m)}( \Lambda_{A_1}(T \oplus Y^{(m)}  ) \Lambda_{A_2}(T \oplus Y^{(m)}  )),
\end{multline}
$a_0(A) = a_1+...+a_s$.
In view of \eqref{Le2-3}, \eqref{Le2-5} and Lemma C, we have
\begin{multline} \nonumber
 \varpi=  \frac{1}{2^{sm}} \sum_{\substack{ t_{i,j} \in \{0,1\} \\ 1 \leq i \leq s,\; 1 \leq j \leq m}}
    \sum_{\substack{  \fm^{[1]}, \fm^{[2]} \in \cP^{\bot,*}_m \\
		\; \rho(\fm^{(\nu)}_i) \leq a_{\nu,i}, \;i=1,...,s, \nu=1,2   }} e\Big(\sum_{i=1}^s \big(\sum_{j=1}^{m} (\fm^{[1]}_{i,j}+
    \fm^{[2]}_{i,j})(
     t_{i,j}+y_{i,j})      \\
          \fm^{[1]}_{i,a_{1,i}}+
    \fm^{[2]}_{i,a_{2,i}} \big) \Big)       =   \sum_{\substack{  \fm^{[1]}, \fm^{[2]} \in \cP^{\bot,*}_m \\
		\; \rho(\fm^{(\nu)}_i) \leq a_{\nu,i}, \;i=1,...,s, \nu=1,2   }}
		\prod_{i=1}^s \prod_{j=1}^m \delta_2(\fm^{[1]}_{i,j}+\fm^{[2]}_{i,j})
      e(\fm^{[1]}_{i,a_{1,i}}+
    \fm^{[2]}_{i,a_{2,i}} ) \\
  =   \sum_{\substack{  \fm^{[1]}, \fm^{[2]} \in \cP^{\bot,*}_m \\
		\; \rho(\fm^{(\nu)}_i) \leq a_{\nu,i}, \;i=1,...,s, \nu=1,2   }}
		\prod_{i=1}^s \prod_{j=1}^m \delta_2(\fm^{[1]}_{i,j}+\fm^{[2]}_{i,j})
		\leq
	\sum_{\substack{  \fm \in \cP^{\bot,*}_m \\
		 \rho(\fm_i) \leq a_{1,i}, \;i=1,...,s   }}	1   = 	\# G_{A_1}.
\end{multline}
Taking into account that $\tau$ does not depend on $Y$, we get from Lemma C and \eqref{Le3-5} that
$A_1=A_2$ and
\begin{multline} \nonumber
  E^{(m)}_{Y,T} ( R^2(T,Y) )= \sum_{A \in I_{m,s,m+1}}
	2^{2m-2s -2a_{1}-\cdots -2a_{s}} \# G_A \\
 \leq
\sum_{\substack{A \in I_m^s \\ a_0(A) \geq m+V_0}}
	2^{2m-2s -2a_{1}-\cdots -2a_{s}} 2^{a_1+ \cdots +a_s - m+\mt_m +1 }=
\sum_{\substack{A \in I_m^s \\ a_0(A) \geq m+V}}
	2^{m-2s -a_{1}-\cdots -a_{s}+\mt_m+ 1} \\
		\leq 2^{-V_0 +\mt_m+1}\sum_{A \in I_m^s} 1 = 2^{-V_0 +\mt_m +1} m^{s} \leq 2^{-10 \log m +1/10\log m} m^s \ll 1/m.
\end{multline}
According to Lemma A, we have
\begin{equation}  \nonumber
       \EE_{Y,T,m} \big( |\cE^{(m)}(\cP_{m}\oplus T , Y ) |^2 \big)      \leq (s 2^{\mt_m})^2
       \leq s^2 2^{1/5 \; \log_2 m} \ll m^{1/5}.
\end{equation}
For $ \mt_m \leq 1/10 \;  \log_2 m $, we obtain
\begin{multline}  \nonumber
\EE_{Y,T,m}((\cE^{(m)}(\cP_{m}\oplus T , Y ) +R(T,Y))^2) \\
\leq 2 \EE_{Y,T,m}((\cE^{(m)}(\cP_{m}\oplus T , Y ))^2) +2 \EE_{Y,T,m}(R(T,Y))^2 \ll  1/m +
    m^{1/5}  \ll
m^{1/5}.
\end{multline}
Therefore Lemma 2 is proved. \qed  \\  \\ \\
%
%
%
{\bf 4. Fourth moments estimates.}\\ \\
{\bf Lemma 3.}  {\it Let   $\mt_m \leq  1/10 \;  \log_2 m$}. Then
\begin{equation}  \nonumber
        \sum_{k=1}^4  E^{(m)}_{Y,T} ( \DD_k^4 )   \ll  m^{2(s-1) -2/5}.
\end{equation} \\
{\bf Proof.}
  Applying Lemma 1, we put
\begin{equation}  \nonumber
  \Psi_A =  r_A(Y) \hat{\Psi}_A, \quad \with \quad \hat{\Psi}_A=  (-1)^{\kappa(A)}  2^{m-s -a_1-\cdots -a_s}    \Lambda_A(T \oplus Y^{(m)}  ).
\end{equation}
By Lemma B and Lemma C, we obtain
\begin{multline} \label{Le3-3}
|\hat{\Psi}_A|=   2^{m-s -a_1-\cdots -a_s} \Bigg|  \sum_{\substack{  \fm \in \cP^{\bot,*}_m, \; \rho(\fm_i) \leq a_i \\ i=1,...,s   }} e\Big(\sum_{i=1}^s \big(\sum_{j=1}^{m} \fm_{i,j}(t_{i,j}+y_{i,j} )
  + \fm_{i,a_{i}} \big) \Big)  \Bigg|  \\
 \leq    2^{m-s -a_1-\cdots -a_s} \#  \{ \fm \in \cP_m^{\bot} \; : \; \rho(\fm_i) \leq a_i, \; i=1,...,s  \} \\
   \leq 2^{m-s -a_1-\cdots -a_s} \; 2^{a_1+ \cdots +a_s - m+\mt_m +1 } = 2^{\mt_m +1 }, \qquad
   |\Psi_A| \leq 2^{\mt_m +1 }.
\end{multline}
From Lemma 1, we get
\begin{multline}\nonumber
     E^{(m)}_{Y,T}  (\DD_k^4(T,Y) )
     = \sum_{A_1,A_2,A_3,A_4 \in I_{m,s,k}}    E^{(m)}_{Y,T} ( \Psi_{A_1} \cdots \Psi_{A_4}) \\
=  \sum_{A_j,A_2,A_3,A_4 \in I_{m,s,k}}  E_Y^{(m)}\big(\prod_{j=1}^{4} r_{A_j}(Y) \sigma \big) , \quad \with \quad  \sigma =
 E_T^{(m)}( \hat{\Psi}_{A_1} \cdots \hat{\Psi}_{A_4}).
\end{multline}
It is easy to see that $\sigma$ does not depend on $Y$. Therefore
\begin{multline}  \nonumber
     E^{(m)}_{Y,T} ( \DD_k^4(T,Y) )  \leq 2^{4(\mt_m +1)}
       \sum_{\substack{  A_j\in I_{m,s,k}, \;  \max_i a_{j,i}=k \\
			 a_0(A_j) \in ( m-\mt_m,m+V_0), j=1,...,4, i=1,...,s}}
			E_Y^{(m)}(\prod_{j=1}^{4} r_{A_j}(Y) )  \\
  \ll (\mt_m +V_0)^4 2^{4(\mt_m +1)}
       \sum_{\substack{  a_{j,i} \in I_{m}, \; \max_i a_{j,i}=k\\ j=1,...,4, \; i=1,...,s-1}}
			E_Y^{(m)}( \prod_{i=1}^{s-1}\prod_{j=1}^{4} r_{a_{j,i}}(Y) ) \\
	\ll (\mt_m +V_0)^4 2^{4 \mt_m }
       \sum_{\substack{  a_{j,i} \in I_{m} \\ j=1,...,4, \; i=1,...,s-2}} \; E_Y^{(m)}( \prod_{i=1}^{s-2}\prod_{j=1}^{4} r_{a_{j,i}}(Y) ) \\
  \ll V_0^4 m^{2/5} \varsigma^{s-2}, \quad \with \quad  \varsigma=
			E_Y^{(m)}( \sum_{  a_{j} \in I_{m}, j=1,...,4} \prod_{j=1}^{4} r_{a_{j,i}}(Y) ).
\end{multline}
Applying Khintchin's inequality (see, e.g., [Skr3, Lemma 3.2]), we get
\begin{equation}  \nonumber
     \varsigma \ll m^2  .
\end{equation}
Hence
\begin{multline}  \nonumber 
    E^{(m)}_{Y,T}  (\DD_k^4(T,Y) ) \ll \log_2^4 m \;  m^{2/5} \; m^{2(s-2)} \\
 \ad \qquad 	  \sum_{k=1}^m   E^{(m)}_{Y,T}  (\DD_k^4(T,Y) )  \ll \log_2^4 m \;  m^{2(s-1) -3/5}	.
\end{multline}		
Therefore Lemma 3 is proved. \qed  \\  \\
Denote by $F_k$ the sigma field in $[0,1)^{2s}$ generated by
\begin{equation}  \nonumber
 \Big\{ \Big[ \frac{\kk_1}{2^k}, \frac{\kk_1 +1}{2^k}\Big)  \times \cdots \times
 \Big[ \frac{\kk_{2s}}{2^k}, \frac{\kk_{2s} +1}{2^k}\Big) \; : \;
  \kk_i =0,1,...,2^k-1, \; i=1,...,2s  \Big\}.
\end{equation}
Denote by $\cF_k$ the sigma field in $[0,1)^{s}$  generated by
\begin{equation}  \nonumber 
 \Big\{ \Big[ \frac{\kk_1}{2^k}, \frac{\kk_1 +1}{2^k}\Big)  \times \cdots \times
 \Big[ \frac{\kk_{s}}{2^k}, \frac{\kk_{s} +1}{2^k}\Big) \; : \;
  \kk_i =0,1,...,2^k-1, \; i=1,...,s  \Big\}.
\end{equation}
We put
\begin{equation} \nonumber
E^{(m)}_{Y,T} ( f(T,Y) \; | \; F_k) = \frac{1}{2^{2s(m-k)}} \sum_{\substack{ t_{i,j}, y_{i,j} \in \{0,1\} \\ 1 \leq i \leq s,\; k+1 \leq j \leq m}} f(T^{(m)}, Y^{(m)})
\end{equation}
and
\begin{equation} \nonumber
E^{(m)}_{Y} ( f(T,Y) \; | \; \cF_k) = \frac{1}{2^{s(m-k)}} \sum_{\substack{  y_{i,j} \in \{0,1\} \\ 1 \leq i \leq s,\; k+1 \leq j \leq m}} f(T^{(m)}, Y^{(m)}),
\end{equation}
\begin{equation} \label{Lem3-5}
E^{(m)}_{T} ( f(T,Y) \; | \; \cF_k) = \frac{1}{2^{s(m-k)}} \sum_{\substack{t_{i,j} \in \{0,1\} \\ 1 \leq i \leq s,\; k+1 \leq j \leq m}} f(T^{(m)}, Y^{(m)}).
\end{equation}
Hence
\begin{equation} \nonumber 
E^{(m)}_{Y,T} ( f(T,Y) \; | \; F_k) = E^{(m)}_{T} \Big( \big(E^{(m)}_{Y} ( f(T,Y) \; \big| \; \cF_k)\big) \;\Big| \;\cF_k) \Big).
\end{equation}
For $A \in I_{m,s,k}$, we have from \eqref{Le3-5} that
$\max_{1 \leq i \leq s}  a_i =k$ and  $ E^{(m)}_{Y} ( r_{A} \; | \; \cF_{k-1}) =0$.
From Lemma 1, we get
\begin{multline}  \nonumber
E^{(m)}_{T} (\Lambda_A(T \oplus Y^{(m)}  ) \; | \; \cF_{k-1} ) =0 \quad \for \quad \exists i\in [1,s] \;\with \; \fm_{i,k} \neq 0, \\
E^{(m)}_{Y} ( r_{A} \Lambda_A(T \oplus Y^{(m)}  ) \; | \; \cF_{k-1} ) =0 \quad \with \quad  \fm_{i,k} = 0 \;\for\; i=1,...,s.
\end{multline}
Now by Lemma 1, we obtain that $\DD_k$ is $F_k$ measurable and that
\begin{multline}   \label{Lem3-7}
E^{(m)}_{Y,T} ( \DD_k \; | \; F_{k-1}) =  \sum_{A \in I_{m,s,k}}(-1)^{\kappa(A)}  2^{m-s -a_1-\cdots -a_s}  \\
\times  E^{(m)}_{Y} \Big(r_A(Y)
\big( E^{(m)}_{T} (\Lambda_A(T \oplus Y^{(m)}  ) \; | \; \cF_{k-1} ) \big) \; | \; \cF_{k-1}  \Big)=0.
\end{multline}
Therefore $(\DD_k, F_k)_{k \geq 1}$ is the martingale difference sequence.\\
Let
\begin{equation}  \label{Lem3-8}
\Omega_m = \sum_{k=1}^m \omega_k \quad \with \quad
\omega_k = E^{(m)}_{Y,T} ( \DD^2_k \; | \; F_{k-1}) -  E^{(m)}_{Y,T} ( \DD^2_k ), \quad \ad \quad
\DD_k = \sum_{A \in I_{m,s,k}} \Psi_A .
\end{equation}
By the martingale CLT (see, e.g., Theorem A), in order to prove Theorem~1 it is enough to verify that Levi's conditional expectation $\Omega_m$ satisfies the bound~:
\begin{equation} \label{Lem3-10}
       E^{(m)}_{Y,T}  (\Omega_m^2)  =o( m^{2(s-1)}).
\end{equation}
 We will prove \eqref{Lem3-10} in Lemma 4 - Lemma 6. \\

By  \eqref{Lem3-8}, we have
\begin{equation}\label{Lem4-3}
\omega_k = \sum_{A_1,A_2 \in I_{m,s,k}} \Theta_{A_1,A_2,k} \;\; \with \;\; \Theta_{A_1,A_2,k}=
E^{(m)}_{Y,T} ( \Psi_{A_1} \Psi_{A_2} \; | \; F_{k-1}) -  E^{(m)}_{Y,T} ( \Psi_{A_1} \Psi_{A_2}  )  .
\end{equation}
 We put
\begin{multline} \label{Lem4-5}
 \cJ_{A_1,A_2,1}= \{i \in \{1,...,s\} \; : \; a_{1,i} =a_{2,i}\}, \quad
         \cJ_{A_1,A_2,2}= \{i \in \{1,...,s\} \; : \; a_{1,i} < a_{2,i}\}, \\
          \cJ_{A_1,A_2,3}= \{i \in \{1,...,s\}  :  a_{1,i}  > a_{2,i}\},\;
d_j = \#  \cJ_{A_1,A_2,j}, \; \quad d_1+d_2+d_3=s.
\end{multline}
It is easy to verify that
\begin{equation}\nonumber
\Theta_{A_1,A_2,k} = \sum_{\substack{d_1+d_2+d_3=s \\ d_i \geq 0, i=1,2,3}} \;\; \sum_{\substack{
J_i \subseteq \{1,...,s\}, \; \#J_i=d_i\\ J_1 \cup J_2 \cup J_3 =\{1,...,s\}, J_i \cap J_j =\emptyset, i\neq j}}   \Theta_{A_1,A_2,k} \prod_{i=1}^3 \d1( \cJ_{A_1,A_2,i} =J_i).
\end{equation}
From \eqref{Lem3-8} and \eqref{Lem4-3},  we have
\begin{equation} \nonumber
\Omega_m = \sum_{\substack{d_1+d_2+d_3=s \\ d_i \geq 0, i=1,2,3}} \;\;  \sum_{\substack{
J_i \subseteq \{1,...,s\}, \; \#J_i=d_i\\ J_1 \cup J_2 \cup J_3 =\{1,...,s\}\\ J_i \cap J_j =\emptyset, \; i\neq j}}\;\;
\sum_{k=1}^m \sum_{A_1,A_2 \in I_{m,s,k}} \Theta_{A_1,A_2,k} \; \prod_{i=1}^3 \d1( \cJ_{A_1,A_2,i} =J_i).
\end{equation}
Applying the Cauchy-Schwarz inequality, we obtain
\begin{multline} \label{Lem4-7}
  E^{(m)}_{Y,T}   (\Omega_m^2 ) \leq s^3 2^{3s} \sum_{\substack{d_1+d_2+d_3=s \\ d_i \geq 0, i=1,2,3}} \;\;  \sum_{\substack{
J_i \subseteq \{1,...,s\}, \; \#J_i=d_i\\ J_1 \cup J_2 \cup J_3 =\{1,...,s\}\\ J_i \cap J_j =\emptyset, \; i\neq j}} \;\; \sum_{k_1, k_2=1}^m  \Upsilon_{d_1,d_2,d_3,J,k}                  \quad \with\\
\Upsilon_{d_1,d_2,d_3,J,k}:=   \sum_{\substack{ A_1,A_2 \in I_{m,s,k_1} \\ A_3,A_4 \in I_{m,s,k_2}}}  E^{(m)}_{Y,T}  (\Theta_{A_1,A_2,k_1}  \Theta_{A_3,A_4,k_2}) \\
   \times \prod_{i=1}^3 \d1( \cJ_{A_1,A_2,i} =\cJ_{A_3,A_4,i} =J_i).
\end{multline}
Hence in order to prove \eqref{Lem3-10} it is enough to verify that
\begin{equation} \label{Lem3-10a}
  \sum_{k_1, k_2=1}^m     \Upsilon_{d_1,d_2,d_3,J,k}  =o( m^{2(s-1)}).
\end{equation} \\ \\
{\bf Lemma 4.}  {\it Let  $d_1=m$ and $\mt_m \leq 1/10 \;  \log_2 m$. Then}
\begin{equation} \label{Lem5-0}
\sum_{k_1, k_2=1}^m   \Upsilon_{m,0,0,J,k} \ll m^{2(s-1)-2/5}.
\end{equation} \\
{\bf Proof.} We consider the case $k_1 \leq k_2$. The proof for the case $k_1 > k_2$  is similar.
By \eqref{Le3-3}, \eqref{Lem4-5} and \eqref{Lem4-7}, we get
\begin{equation}\label{Lem5-1}
    |\Psi_A| \leq 2^{\mt_m+1}   \quad \ad \quad A_1=A_2, \quad    A_3=A_4 \quad \for \quad d_1=m.
\end{equation}
Let
\begin{multline} \label{Lem4-80}
 \fg_0:=  \sum_{m \geq k_2 \geq  k_1 \geq 1}   |\Upsilon_{m,0,0,J,k}| = \fg_1 + \fg_2, \;\; \fg_1 = \sum_{m \geq k_2 \geq  k_1 \geq 1} |\Upsilon_{m,0,0,J,k}|  \d1( k_2 -k_1 \leq 10s V_0),\\
  \fg_2= \sum_{m \geq k_2 \geq  k_1 \geq 1}
  |\Upsilon_{m,0,0,J,k}|  \d1( k_2 -k_1 > 10s V_0).
\end{multline}
Let us consider  $\fg_1$. \\
From \eqref{Le3-5}, \eqref{Lem4-3} - \eqref{Lem4-7} and   \eqref{Lem5-1}, we get
\begin{multline}   \label{Lem5-10}
 \Upsilon_{m,0,0,J,k}  \leq 2^{4\mt_m+4} \sum_{\substack{ A_1,A_2 \in I_{m,s,k_1} \\ A_3,A_4 \in I_{m,s,k_2}}}  \d1(A_1=A_2) \; \d1(A_3=A_4)\; \d1( \max_i a_{1,i} =k_1)  \\
  \times \d1( \max_i a_{3,i} =k_2)\; \prod_{j=1}^4 \d1( \sum_{i=1}^s a_{j,i} \in (m-t,m+V_0)) \\
  \leq 2^{2/5 \log_2 m +4} \sum_{\substack{ a_{j,i} \in \{0,...,m\} \\ j=1,2,\; i=1,...,s}}
\prod_{j=1}^2 \d1( \max_i a_{j,i} =k_2) \prod_{j=1}^4 \d1( \sum_{i=1}^s a_{j,i} \in (m-t,m+V_0)) \\
  \leq 16 m^{2/5} m^{2(s-2)} (\mt_m+V_0)^2 \ll
  m^{2(s-2)+2/5} \log_2^2 m, \quad \ad \quad \fg_1 \ll    m^{2(s-1)-3/5} \log_2^3 m.
\end{multline}
Now we will prove that  $\fg_2=0$. \\

According to \eqref{Lem4-7},  it is enough to verify that
\begin{equation}\label{Lem5-11}
      E^{(m)}_{Y,T}  (\Theta_{A_1,A_2,k_1}  \Theta_{A_3,A_4,k_2}) =0.
\end{equation}
We have
\begin{multline}   \nonumber
  E^{(m)}_{Y,T} (\Theta_{A_1,A_2,k_1}  \Theta_{A_3,A_4,k_2}) =  E^{(m)}_{Y,T}  (\Gamma),
      \quad \with \\
       \Gamma =
       E^{(m)}_{Y,T} (\Theta_{A_1,A_2,k_1}  \Theta_{A_3,A_4,k_2}\; |\; F_{k_1-1}).
\end{multline}
By  \eqref{Lem4-3}, we have that $\Theta_{A_1,A_2,k_1} $ is $F_{k_1-1}$ measurable and
\begin{equation}\nonumber
      \Gamma =
       \Theta_{A_1,A_2,k_1} E^{(m)}_{Y,T} (  \Theta_{A_3,A_4,k_2}\; |\; F_{k_1-1}).
\end{equation}
From  \eqref{Lem5-1}, we obtain that $A_3=A_4$ for $d_1=m$.
By \eqref{Lem4-3}, we get that in order to prove \eqref{Lem5-11}, it is enough to verify that
\begin{equation}   \nonumber 
      E^{(m)}_{Y,T} (  \Theta_{A_3,A_4,k_2}\; |\; F_{k_1-1})=
E^{(m)}_{Y,T} ( \Psi^2_{A_3} \; | \; F_{k_1-1}) -  E^{(m)}_{Y,T} ( \Psi^2_{A_3}  ) =0.
\end{equation}
In view of Lemma 1, we have
\begin{multline}\nonumber
  \wp:=      2^{2(-m+s +a_1+\cdots +a_s)}  \Psi^2_{A_3} =   \Lambda^2_{A_3}(T \oplus Y^{(m)}  )\\
= \sum_{\substack{  \fm^{[1]},\fm^{[2]} \in \cP^{\bot,*}_m, \; \rho(\fm^{[1]}_i) \leq a_{3,i} \\\rho(\fm^{[2]}_i) \leq a_{3,i} ,\;  i=1,...,s   }} e\Big(\sum_{i=1}^s  \big(\sum_{j=1}^{m} (\fm^{[1]}_{i,j}+\fm^{[2]}_{i,j})(t_{i,j}+ y_{i,j})
    + \fm^{[1]}_{i,a_{3,i}}+\fm^{[2]}_{i,a_{3,i}}   \big) \Big).
\end{multline}
We want to verify that
\begin{equation}\label{Lem5-13}
E^{(m)}_{Y,T} ( \wp \; | \; F_{k_1-1}) -  E^{(m)}_{Y,T} (\wp  ) =0.
\end{equation}
By \eqref{Le2-3} and \eqref{Lem3-5}, we get
\begin{multline}\nonumber
E^{(m)}_{Y,T} ( \wp \; | \; F_{k_1-1}) = \sum_{\substack{  \fm^{[1]},\fm^{[2]} \in \cP^{\bot,*}_m, \; \rho(\fm^{[1]}_i) \leq a_{3,i} \\\rho(\fm^{[2]}_i) \leq a_{3,i} ,\;  i=1,...,s   }}e\Big(\sum_{i=1}^s  \big(\sum_{j=1}^{m} (\fm^{[1]}_{i,j}+\fm^{[2]}_{i,j})(t_{i,j}+ y_{i,j})
    + \fm^{[1]}_{i,a_{3,i}}\\
  +\fm^{[2]}_{i,a_{3,i}}   \big) \Big)   \prod_{i=1}^s \prod_{i=k_1}^m \d1(\fm^{[1]}_{i,j}=\fm^{[2]}_{i,j}).
\end{multline}
Hence
\begin{multline}\label{Lem5-14}
E^{(m)}_{Y,T} ( \wp \; | \; F_{k_1-1}) - E^{(m)}_{Y,T} ( \wp ) = \sum_{\substack{  \fm^{[1]},\fm^{[2]} \in \cP^{\bot,*}_m, \; \rho(\fm^{[1]}_i) \leq a_{3,i} \\\rho(\fm^{[2]}_i) \leq a_{3,i} ,\;  i=1,...,s   }} e\Big(\sum_{i=1}^s  \big(\sum_{j=1}^{m} (\fm^{[1]}_{i,j}+\fm^{[2]}_{i,j})(t_{i,j}+ y_{i,j}) \\
    + \fm^{[1]}_{i,a_{3,i}}+\fm^{[2]}_{i,a_{3,i}}   \big) \Big) \prod_{i=1}^s \d1\Big(  \rho(\fm^{[1]}_i -\fm^{[2]}_{i}) \leq \min(a_{3,i},k_1)\Big)   \d1(\fm^{[1]} \neq \fm^{[2]}).
\end{multline}
For $A_3 \in I_{m,s,k} $, we obtain from  \eqref{Le3-5}
\begin{equation}\nonumber
a_{3,1}+ \cdots +a_{3,s} \leq m+V_0, \quad \ad \quad \max_{1 \leq i \leq s}(a_{3,i})=k_2, \quad \exists i_0 \; \with \; a_{3,i_0}=k_2.
\end{equation}
For $k_2 -k_1 \geq 10sV_0$, we get
\begin{multline}  \nonumber 
  \sum_{1 \leq i \leq s}  \min(a_{3,i},k_1) =  \sum_{1 \leq i \leq s, i \neq i_0}  \min(a_{3,i},k_1) +  \min(k_2,k_1) \leq  \sum_{1 \leq i \leq s}  a_{3,i}  -(k_2 -k_1) \\
  \leq m+V_0 -(k_2 -k_1) \leq m-V_0, \quad V_0 =10 \; \log_2 m .
\end{multline}
Suppose that $ \fm^{[1]} -\fm^{[2]} \neq 0 $.
Bearing in mind that $ \fm^{[1]} -\fm^{[2]} \in \cP^{\bot}_m  $, we have
\begin{equation}\nonumber
m-V_0 \geq   \sum_{i=1}^s \min(a_{3,i},k_1) \geq \sum_{i=1}^s \rho(\fm^{[1]}_i -\fm^{[2]}_{i}) \geq m-\mt_m \quad \for \quad  \mt_m < 1/10 \;   \log_2 m.
\end{equation}
We  have a contradiction. Therefore the sum in \eqref{Lem5-14} is zero and \eqref{Lem5-11} - \eqref{Lem5-13} are true. Hence $\fg_2 =0$.
By  \eqref{Lem4-80} and \eqref{Lem5-10}, Lemma 4  is proved. \qed \\ \\
Let
\begin{equation} \label{Lem6-1}
     \ff_{A_1,A_2,k_1}  =
       \begin{cases}
    1,  & \; {\rm if}  \; \exists \; i \in J_2 \cup J_3, \; j \in \{1,2\} \;\; \with \;\; a_{j,i} =k_1,\\
    0, &{\rm otherwise}
  \end{cases}
\end{equation}
and let
\begin{equation} \nonumber
     \ff^{'}_{A_3,A_4,k_2}  =
       \begin{cases}
    1,  & \; {\rm if}  \; \exists \; i \in J_2 \cup J_3, \;\; \with \;\; \max(a_{3,i},a_{4,i}) \geq k_1,\\
    0, &{\rm otherwise}
  \end{cases} .
\end{equation}\\ \\
{\bf Lemma 5.}  {\it Let  $d_1<m$ and   $\mt_m \leq 1/10 \; \log_2 m$. Then}
\begin{multline}\label{Lem6-2}
      \Upsilon_{d_1,d_2,d_3,J,k} \ll    \sum_{\substack{ A_1,A_2 \in I_{m,s,k_1} \\ A_3,A_4 \in I_{m,s,k_2}}} \cL_{A,k} \ll m^{2(s-2) +3/5}, \\
      \with \qquad   \cL_{A,k} = \sL_{A,k}\;  W_1 \; W_2 \; W_3 \; W_4 \; W_5, \quad \sL_{A,k} \leq  m^{2/5}, \\
%
\where \qquad \qquad \qquad \qquad \qquad \qquad \qquad \qquad \qquad \qquad \qquad \qquad \qquad \qquad \qquad \qquad \qquad \\
W_1 = \prod_{i \in J_1} \d1(a_{1,i} = a_{2,i})  \d1(a_{3,i} = a_{4,i}), \quad
W_2 = \prod_{i \in J_2} \d1(a_{1,i} < a_{2,i})  \d1(a_{3,i} < a_{4,i}),  \\
%
W_3 = \prod_{i \in J_3} \d1(a_{1,i} > a_{2,i})  \d1(a_{3,i} > a_{4,i}), \quad
W_4 = \prod_{i \in J_2 \cup J_3} \d1(a_{1,i} = a_{3,i})  \d1(a_{2,i} = a_{4,i}), \\
W_5= \d1(\ff_{A_1,A_2,k_1}=0)  \;  \d1(\ff^{'}_{A_3,A_4,k_2}=0), \\
  \sL_{A,k}:=
2^{a_A} \sum_{\substack{ \fm^{[\ell]}  \in \cP^{\bot,*}_m, \; \rho(\fm^{[l]}_i) \leq a_{l,i} \\  \ell=1,...,4,\; i=1,...,s   }}  \d1\big(\sum_{\ell=1}^4  \fm^{[\ell]} =0 \big), \quad
   a_A=4m - 4s - \sum_{j=1}^4 \sum_{i=1}^s  a_{j,i}.
\end{multline}  \\
{\bf Proof.} We consider the case $k_1 \leq k_2$. The proof for the case $k_1 > k_2$  is similar.

Firstly, we will consider $ \Theta_{A_1,A_2,k_1}$. We have
\begin{multline}\label{Lem6-7}
\hat{\Theta}:= E^{(m)}_{Y,T} \big(r_{A_1} (Y)r_{A_2}(Y)  \Lambda_{A_1}(T \oplus Y^{(m)})  \Lambda_{A_2}(T \oplus Y^{(m)})  \; | \; F_{k_1-1}\big)\\
  =E^{(m)}_{Y} \Bigg( r_{A_1}(Y)r_{A_2}(Y) \Big(E^{(m)}_{T}\big(\Lambda_{A_1}(T \oplus Y^{(m)}) \Lambda_{A_2}(T \oplus Y^{(m)} ) \; | \; \cF_{k_1-1} \big)\Big) \;|\; \cF_{k_1-1} \Bigg).
\end{multline}
Using Lemma 1 and \eqref{Lem3-5}, we derive
\begin{multline}\label{Lem6-6}
E^{(m)}_{T} \big( \Lambda_{A_1}(T \oplus Y^{(m)})  \Lambda_{A_2}(T \oplus Y^{(m)})   |  \cF_{k_1-1}\big) =
\sum_{\substack{  \fm^{[1]},\fm^{[2]} \in \cP^{\bot,*}_m, \ell=1,2 \\
 \rho(\fm^{[\ell]}_i) \leq a_{\ell,i} ,  i=1,...,s   }} e\Big(\sum_{i=1}^s \big( \sum_{j=1}^{m} (\fm^{[1]}_{i,j} \\
   +\fm^{[1]}_{i,j})  (t_{i,j}+ y_{i,j}) + \fm^{[1]}_{i,a_{1,i}} + \fm^{[2]}_{i,a_{2,i}} \big) \Big) \prod_{i=1}^s \prod_{i=k_1}^m  \d1(\fm^{[1]}_{i,j}= \fm^{[2]}_{i,j}).
\end{multline}
It is easy to see that the left side of this equality does not depends on $y_{i,j}$ for $i=1,...,s$ and $j \geq k_1$.
By \eqref{Lem4-5}, we have $A_1, A_2 \in I_{m,s,k_1}$ and $ \max_{j=1,2,\; i=1,...,s} a_{j,i} =k_1$.
For  $\ff_{A_1,A_2,k_1}=1$ and $d_1 <m$, we get from \eqref{Lem4-5} and  \eqref{Lem6-1}  that there exists $i_0 \in \cJ_2 \cJ_3$ such that $a_{1,i_0} \neq a_{2,i_0} $
 and  $ \max_{j=1,2} a_{j,i_0}=k_1$.
Therefore
\begin{equation}  \nonumber 
  \sum_{y_{i_0,k_1} =0,1} r_{a_{1,i_0}} (y_{i_0,k_1})  r_{a_{2,i_0}} (y_{i_0,k_1})=0 \quad \ad \quad
    \hat{\Theta}=0.
\end{equation}
Similarly, we derive from \eqref{Lem6-7} and definition  \eqref{Lem4-5} of $d_1$ that
\begin{equation}  \nonumber 
E^{(m)}_{Y,T} \big(r_{A_1} (Y)r_{A_2}(Y)  \Lambda_{A_1}(T \oplus Y^{(m)})  \Lambda_{A_2}(T \oplus Y^{(m)})  \big) =0 \quad \for \quad d_1 <m.
\end{equation}
Substituting this equality and the equality $\hat{\Theta}=0 $ in \eqref{Lem4-3},  we obtain
\begin{equation}   \nonumber 
 \Theta_{A_1,A_2,k_1}=
E^{(m)}_{Y,T} ( \Psi_{A_1} \Psi_{A_2} \; | \; F_{k_1-1})
= \d1(\ff_{A_1,A_2,k_1}=0 ) E^{(m)}_{Y,T} ( \Psi_{A_1} \Psi_{A_2}  \; | \; F_{k_1-1}).
\end{equation}
We put
\begin{equation}\label{Lem6-10}
\vartheta_{A_1,A_2,k_1}
=  E^{(m)}_{T} ( \tilde{\Psi}_{A_1} \tilde{\Psi}_{A_2}  \; | \; \cF_{k_1-1}), \quad \tilde{\Psi}_{A} = r_{A} \Psi_{A} .
\end{equation}
 For $\ff_{A_1,A_2,k_1}=0$, we have from \eqref{Lem6-1} and \eqref{Lem6-6} that $ \vartheta_{A_1,A_2,k_1}$ and $r_{A_1} (Y)r_{A_2}(Y)   $
does not depend on $y_{i,j}$  ($i=1,...,s$, $j \geq k_1$).

Hence
\begin{multline}\label{Lem6-11}
\Theta_{A_1,A_2,k_1}= \d1(\ff_{A_1,A_2,k_1}=0 )  E^{(m)}_{Y} ( r_{A_1} (Y)r_{A_2}(Y) \vartheta_{A_1,A_2,k_1} \; | \; \cF_{k_1-1})\\
  =\d1(\ff_{A_1,A_2,k_1}=0 ) r_{A_1}(Y)r_{A_2}(Y)   \vartheta_{A_1,A_2,k_1}.
\end{multline}
Let
\begin{equation}     \nonumber 
  \check{\cL}_{A,k} :=     E^{(m)}_{Y,T}  \big(\Theta_{A_1,A_2,k_1}  \Theta_{A_3,A_4,k_2} \big).
\end{equation}
From \eqref{Lem4-5}, we obtain
\begin{equation}\label{Lem6-11b}
          \prod_{i=1}^3 \d1( \cJ_{A_1,A_2,i} =\cJ_{A_3,A_4,i} =J_i) =
          W_1 \; W_2 \; W_3 .
\end{equation}
By \eqref{Lem4-7}, we have
\begin{multline} \label{Lem6-11f}
\Upsilon_{d_1,d_2,d_3,J,k} =   \sum_{\substack{ A_1,A_2 \in I_{m,s,k_1} \\ A_3,A_4 \in I_{m,s,k_2}}}
   \hat{\cL}_{A,k}, \;\; \with \;\; \hat{\cL}_{A,k}:=  \check{\cL}_{A,k}  \prod_{i=1}^3 \d1( \cJ_{A_1,A_2,i} =\cJ_{A_3,A_4,i} =J_i)\\
 =
   \check{\cL}_{A,k}  W_1 \; W_2 \; W_3.
\end{multline}
Taking into account that $\Theta_{A_1,A_2,k_1} $ is $F_{k_1-1}$ measurable (see \eqref{Lem4-3}), we get
\begin{equation}\label{Lem6-11c}
  \check{\cL}_{A,k} =
      E^{(m)}_{Y,T} \big(\Theta_{A_1,A_2,k_1}  \tilde{\Theta}_{A_3,A_4,k_2} \big) \;\; \with \;\;  \tilde{\Theta}_{A_3,A_4,k_2} =  E^{(m)}_{Y,T}  \big( \Theta_{A_3,A_4,k_2} \;|\; F_{k_1-1} \big).
\end{equation}
Using the equality $ E^{(m)}_{Y,T} (E^{(m)}_{Y,T} ( f  |  F_{k_2-1}) \; | \;  F_{k_1-1}) =
   E^{(m)}_{Y,T} ( f  |  F_{k_1-1})$, we get from \eqref{Lem4-3} that
\begin{equation}    \nonumber 
\tilde{\Theta}_{A_3,A_4,k_2}=
E^{(m)}_{Y,T} ( \Psi_{A_3} \Psi_{A_4} \; | \; F_{k_1-1}) -  E^{(m)}_{Y,T} ( \Psi_{A_3} \Psi_{A_4}  )  .
\end{equation}
Analogously to \eqref{Lem6-6} - \eqref{Lem6-11}, we get
\begin{multline}\label{Lem6-12}
\tilde{\Theta}_{A_3,A_4,k_2} = \d1(\ff^{'}_{A_3,A_4,k_2}=0 ) E^{(m)}_{Y} ( r_{A_3} (Y)r_{A_4}(Y) \tilde{\theta}_{A_3,A_4,k_2} \; | \; \cF_{k_1-1})  \\
  = \d1(\ff^{'}_{A_3,A_4,k_2}=0 ) r_{A_3}(Y)r_{A_4}(Y) \tilde{\theta}_{A_3,A_4,k_2},
\end{multline}
where
\begin{equation}    \nonumber 
\tilde{\theta}_{A_3,A_4,k_2}
:=  E^{(m)}_{T} ( \tilde{\Psi}_{A_3} \tilde{\Psi}_{A_4}  \; | \; \cF_{k_1-1}).
\end{equation}
Substituting \eqref{Lem6-11}  and  \eqref{Lem6-12}  in \eqref{Lem6-11c}, we obtain
\begin{multline}\label{Lem6-13ac}
\check{\cL}_{A,k} =   \d1(\ff_{A_1,A_2,k_1}=0 ) \;
\d1(\ff^{'}_{A_3,A_4,k_2}=0 )   \;E^{(m)}_{Y}  \big( \prod_{j=1}^4 r_{A_j} (Y) \cL_1 \big) \\
 =W_5 E^{(m)}_{Y}  \big( \prod_{j=1}^4 r_{A_j} (Y) \cL_1 \big) ,
 \with \qquad
  \cL_1 := E^{(m)}_{T} ( \vartheta_{A_1,A_2,k_1} \;  \tilde{\theta}_{A_3,A_4,k_2} ).
\end{multline}
By Lemma 1, we get
\begin{equation}   \nonumber 
    \tilde{\Psi}_A =  (-1)^{\kappa(A)}  2^{m-s - a_0(A)}   \Lambda_A(T \oplus Y^{(m)}  ), \quad
    a_0(A) = a_1+ \cdots  +a_s.
\end{equation}
From \eqref{Lem6-10} and \eqref{Lem6-6} , we get
\begin{multline}   \nonumber
\vartheta_{A_1,A_2,k_1}
= (-1)^{\kappa(A_1) + \kappa(A_2)}  2^{2m-2s - a_0(A_1) - a_0(A_2)}
E^{(m)}_{T} \big( \Lambda_{A_1}(T \oplus Y^{(m)}  )\Lambda_{A_2}(T \oplus Y^{(m)}  )  \\ | \; \cF_{k_1-1}\big)
=(-1)^{\kappa(A_1) + \kappa(A_2)}  2^{2m-2s - a_0(A_1) - a_0(A_2)}
\sum_{\substack{  \fm^{[1]},\fm^{[2]} \in \cP^{\bot,*}_m, \; \rho(\fm^{[\ell]}_i) \leq a_{\ell,i} \\ \ell=1,2,\;  i=1,...,s   }} e\Big(\sum_{i=1}^s \big( \sum_{j=1}^{m} (\fm^{[1]}_{i,j} \\
   +\fm^{[1]}_{i,j})  (t_{i,j}+ y_{i,j}) + \fm^{[1]}_{i,a_{1,i}} + \fm^{[2]}_{i,a_{2,i}} \big) \Big) \prod_{i=1}^s \prod_{i=k_1}^m \d1(\fm^{[1]}_{i,j}= \fm^{[2]}_{i,j}).
\end{multline}
We obtain a similar expression for $\tilde{\theta}_{A_3,A_4,k_2} $. Thus, we get from  \eqref{Lem6-13ac} that
\begin{equation}    \nonumber 
 |\cL_1 |  \leq   2^{a_A} \sum_{\substack{ \fm^{[\ell]}  \in \cP^{\bot,*}_m, \; \rho(\fm^{[l]}_i) \leq a_{l,i} \\  \ell=1,...,4,\; i=1,...,s   }}  | \cL_2|  \prod_{i=1}^s \prod_{i=k_1}^m \d1(\fm^{[1]}_{i,j}=\fm^{[2]}_{i,j}) \d1(\fm^{[3]}_{i,j}=\fm^{[4]}_{i,j})
\end{equation}
where
\begin{multline}\label{Lem6-14}
\cL_2=E^{(m)}_{T} \Big(   e\Big(\sum_{\ell=1}^4 \sum_{i=1}^s \big( \sum_{j=1}^{m}
   \fm^{[\ell]}_{i,j} (t_{i,j}+ y_{i,j}) +  \fm^{[\ell]}_{i,a_{\ell,i}} \big) \Big) \Big), \\
   a_A=4m - 4s - \sum_{j=1}^4 \sum_{i=1}^s  a_{j,i}.
\end{multline}
Let us consider the case $\check{\cL}_{A,k}  \neq 0$.
From \eqref{Lem6-13ac} and \eqref{Lem6-14}, we obtain
\begin{equation} \label{Lem6-13a1}
E^{(m)}_{Y}  \big( \prod_{j=1}^4 r_{A_j} (Y)\big) =1 \quad \ad \quad
\sum_{\ell=1}^4  \fm^{[\ell]} =0.
\end{equation}
Hence, we get from \eqref{Lem4-5} that
\begin{multline}   \nonumber 
a_{1,i} = a_{2,i}, \;\;a_{3,i} = a_{4,i} \;\; \for \; i \in J_1, \quad
a_{1,i} = a_{3,i}, \;\;\;a_{2,i} = a_{4,i} \;\; \for \; i \in J_2 \cup J_3, \; W_4=1,\\
a_{1,i} < a_{2,i}, \quad a_{3,i} < a_{4,i} \quad \for \;\; i \in J_2, \qquad  a_{1,i} > a_{2,i},\quad a_{3,i} > a_{4,i} \quad \for \;\; i \in J_3.
\end{multline}
Applying \eqref{Lem6-11f} and   \eqref{Lem6-13ac}, we get
\begin{equation}\label{Lem6-13a}
 |\hat{\cL}_{A,k}| \leq W_1 W_2 \cdots W_5 \; |\cL_1 |.
\end{equation}
Using \eqref{Lem6-2} and \eqref{Lem6-13a1}, we have
\begin{equation}\label{Lem6-14ab}
 |\cL_1| \leq  \sL_{A,k} = 2^{a_A} \sum_{\substack{ \fm^{[\ell]}  \in \cP^{\bot,*}_m, \; \rho(\fm^{[l]}_i) \leq a_{l,i} \\  \ell=1,...,4,\; i=1,...,s   }}  \d1\big(\sum_{\ell=1}^4  \fm^{[\ell]} =0 \big).
\end{equation}
In view of Lemma B and Lemma C, we get
\begin{multline}\label{Lem6-14a}
     \sL_{A,k} \leq     2^{a_B}, \;\with \; a_B = a_A +\sum_{j=1}^4 (a_0(A_j) -\brho(\PP_m^{\bot})+1)= 4(m-s)-  \sum_{j=1}^4 \brho(\PP_m^{\bot}) \\
   \leq 4(\mt_m-s+1) \leq 2/5 \; \log_2 m.
\end{multline}
By \eqref{Le3-6}, we have
\begin{equation}    \nonumber
  \# I_{m,s,k} \leq m^{s-2}(t+V_0).
\end{equation} \
From   \eqref{Lem6-11f}, \eqref{Lem6-13a}, \eqref{Lem6-14ab} and   \eqref{Lem6-14a}, we obtain
\begin{equation} \nonumber
  \sL_{A,k} \leq  m^{2/5} \quad \ad \qquad    \Upsilon_{d_1,d_2,d_3,J,k} \ll    \sum_{\substack{ A_1,A_2 \in I_{m,s,k_1} \\ A_3,A_4 \in I_{m,s,k_2}}} \sL_{A,k}\;  W_1  W_2 \cdots W_5 \ll m^{2(s-2) +3/5}.
\end{equation}
Hence Lemma 5  is proved. \qed \\ \\
{{\bf Lemma 6.}  {\it Let $d_1 <m$,  $\mt_m \leq  1/10 \; \log_2 m$}. Then
\begin{equation}  \nonumber 
   \sum_{k_1, k_2=1}^m    \Upsilon_{d_1,d_2,d_3,J,k} \ll m^{2(s-1)-1/5}.
\end{equation} \\
{\bf Proof.} We consider the case $k_2 \geq  k_2$. The proof for the case $k_1 > k_2$  is similar.
Firstly, we will consider the case $0 \leq k_2 -k_1 \leq 10s V_0, \;V_0 =10 \log_2 m $.

Let
\begin{multline} \label{Lem7-80}
 \fg_3=  \sum_{m \geq k_2 \geq  k_1 \geq 1}   \Upsilon_{d_1,d_2,d_3,J,k} = \fg_4 + \fg_5, \;\; \fg_4 = \sum_{m \geq k_2 \geq  k_1 \geq 1} |\Upsilon_{d_1,d_2,d_3,J,k}|  \d1( k_2 -k_1 \leq 10s V_0),\\
  \fg_5= \sum_{m \geq k_2 \geq  k_1 \geq 1} \fg_{5,k_1,k_2} \quad \with \quad
\fg_{5,k_1,k_2}:=   |\Upsilon_{d_1,d_2,d_3,J,k}|  \d1( k_2 -k_1 > 10s V_0).
\end{multline}

Applying  Lemma 5,   we get
\begin{equation} \label{Lem7-85}
\fg_4  \ll
  \sum_{\substack{ m \geq k_2 \geq  k_1 \geq 1\\ k_2  \leq k_1 + 2V_0}}  m^{2(s-2) +3/5 }  \ll m^{2(s-1) -2/5} \log_2 m.
\end{equation}

Let us consider the case $k_2 -k_1 \geq 10s V_0$.
We will use Lemma 5 :

\begin{equation}  \nonumber 
     \Upsilon_{d_1,d_2,d_3,J,k} \ll    \sum_{\substack{ A_1,A_2 \in I_{m,s,k_1} \\ A_3,A_4 \in I_{m,s,k_2}}}  \sL_{A,k}\;  W_1  W_2\cdots  W_5.
\end{equation}

We put $\bA =(A_1,A_2,A_3,A_4)$,
\begin{equation}  \nonumber 
    J_4=  \cJ_{4,\bA}= \{i \in  \cJ_{A_1,A_2,1} \; : \; a_{1,i}  > a_{3,i}\}, \quad
         J_5=\cJ_{5,\bA}= \{i \in  \cJ_{A_1,A_2,1} \; : \; a_{1,i}  \leq a_{3,i}\}
\end{equation}
and
\begin{equation} \label{Lem7-1}
\fa_{j,\nu} =\sum_{i \in \cJ_{\nu}}
a_{j,i}, \quad \;\; \fa_{j,\nu}=0 \;\; \for \;\; \cJ_{\nu}=\emptyset,  \quad \quad
d_4 = \# \cJ_{4},\;\;d_5 =\# \cJ_{5}.
\end{equation}
From Lemma 5, \eqref{Le3-5}, \eqref{Lem4-5} and \eqref{Lem6-1} for $k_2 -k_1 > 10s V_0$, we have that if \\$ W_1 \; W_2 \cdots W_5 =1 $, then
\begin{multline} \label{Lem7-2}
a_{1,i} = a_{2,i}, \;\;\;a_{3,i} = a_{4,i} \;\;\;\;\; \for \;\; i \in J_1, \quad \quad
a_{1,i} = a_{3,i}, \;\;\;a_{2,i} = a_{4,i} \;\; \;\; \for \;\; i \in J_2 \cup J_3,\\
a_{1,i} < a_{2,i} \quad \for \;\; i \in J_2, \qquad \quad a_{1,i} > a_{2,i} \quad \for \;\; i \in J_3,   \\
a_{1,i} > a_{3,i} \quad\for \;\; i \in \cJ_{4,\bA}, \qquad \quad
a_{1,i} \leq a_{3,i} \quad\for \;\; i \in \cJ_{5,\bA}, \\
  \max_{j=1,2, i \in [1,s]} a_{j,i}=k_1, \quad  \max_{j=3,4, i \in [1,s]} a_{j,i}=k_2, \quad
  \max_{i \in J_2 \cup J_3, j=1,2} a_{j,i} =   \max_{i \in J_2 \cup J_3, j=3,4} a_{j,i} , \\
\max_{i \in J_1} a_{1,i} \leq k_1, \quad   \quad  \ff^{'}_{A_3,A_4,k_2}=0,\quad  \max_{i \in J_1} a_{3,i} =k_2, \quad
\fa_{3,5} - \fa_{1,5} \geq k_2-k_1,  \\
  \cJ_{4,\bA} \cup  \cJ_{5,\bA} = J_{1}, \quad  a_{j,1}+ \cdots +  a_{j,s}   \in ( m-\mt_m,m+V_0) \quad
\for \quad A_j \in I_{m,s,k}.
\end{multline}
Hence  $\fa_{1,4}=\fa_{2,4}$, $\fa_{1,5}=\fa_{2,5}$,
\begin{equation} \label{Lem7-2a}
\fa_{1,1} = \fa_{2,1}, \;\; \fa_{3,1} = \fa_{4,1}, \;\;
        \fa_{1,2} = \fa_{3,2}, \;\;  \fa_{1,3} = \fa_{3,3}, \;\;  \fa_{2,2} = \fa_{4,2}, \;\;  \fa_{2,3} = \fa_{4,3}.
\end{equation}
We get from \eqref{Lem7-1} and \eqref{Lem7-2} that
\begin{equation} \label{Lem7-2b}
\sum_{\ell=2}^5 \fa_{j,\ell} \in (m-\mt_m, m+V_0), \quad \;\;\;
\sum_{\ell=2}^5 \fa_{j_1,\ell} -  \sum_{\ell=2}^5 \fa_{j_2,\ell} \in (-\mt_m -V_0, \mt_m + V_0).
\end{equation}
Taking into account that $\fa_{1,4}=\fa_{2,4}$, $\fa_{1,5}=\fa_{2,5}$, we get from \eqref{Lem7-2b}
\begin{equation}  \label{Lem7-30}
(\fa_{1,2} + \fa_{1,3}) - (\fa_{2,2} + \fa_{2,3})= (\fa_{1,3} - \fa_{2,3}) - (\fa_{2,2} - \fa_{1,2}) = \hat{\fc}_{1} -  \fc_{1} \in (-\mt_m -V_0, \mt_m + V_0) ,
\end{equation}
where
\begin{equation} \label{Lem7-32}
\hat{\fc}_{1} = \fa_{1,3} - \fa_{2,3}, \quad \ad \quad  \fc_{1} = \fa_{2,2} - \fa_{1,2}.
\end{equation}
Taking into account that $\fa_{1,2}=\fa_{3,2}$, $\fa_{1,3}=\fa_{3,3}$, we get from \eqref{Lem7-2b}
\begin{equation} \label{Lem7-34}
(\fa_{3,4} + \fa_{3,5}) - (\fa_{1,4} + \fa_{1,5})= (\fa_{3,5} - \fa_{1,5})  -(\fa_{1,4} - \fa_{3,4}) =\hat{\fc}_{2} - \fc_{2}                \in (-\mt_m -V_0, \mt_m + V_0),
\end{equation}
where
\begin{equation} \label{Lem7-36}
         \hat{\fc}_{2}  = \fa_{3,5} - \fa_{1,5},  \quad  \fc_{2} = \fa_{1,4} - \fa_{3,4}.
\end{equation}
In view of Lemma B, Lemma C, Lemma 5, \eqref{Lem7-1} and  \eqref{Lem7-2}, we get for $\fm^{[j]} \in   G_{A_j}$ that
\begin{multline} \nonumber
 m+V_0 \geq \sum_{i=1}^s \fa_{j,i} , \;\; \sum_{i=1}^s \rho(\fm^{[j]}_{i}) \geq m -\mt_m, \;\;
 \fa_{j,i} - \rho( \fm^{[j]}_{i})  \geq 0, \quad
\sum_{i=1}^s (\fa_{j,i} - \rho( \fm^{[j]}_{i})) \geq 0, \\
 -m +\mt_m \geq  -\sum_{i=1}^s \rho(\fm^{[j]}_{i}) , \;\; V_0 +\mt_m \geq \sum_{i=1}^s (\fa_{j,i} - \rho( \fm^{[j]}_{i})) \geq 0, \; \;   V_0 +\mt_m  \geq \fa_{j,i} - \rho( \fm^{[j]}_{i}) \geq 0.
\end{multline}
Hence, there exist $\epsilon_{j,i} \in [0,1]$ with
\begin{equation} \label{Lem7-2ab}
      \fa_{j,i} =  \rho( \fm^{[j]}_{i}) + \epsilon_{j,i} (V_0 +\mt_m ), \quad i=1,...,s,\; j=1,...,4.
\end{equation}\\

  Let us consider the case $d_1=0$.\\
 Bearing in mind  that $A_1,A_2 \in I_{m,s,k_1}$, we get  $\max_i(a_{1,i}) =k_1$.
Using Lemma 5, we obtain  $\ff_{A_1,A_2,k_1} =0$. By \eqref{Lem6-1}, $\max_i(a_{1,i}) <k_1$  for $d_1=0$. We have a contradiction. Therefore
\begin{equation}  \nonumber 
     \Upsilon_{0,d_2,d_3,J,k}  =0.
\end{equation}

 Let us consider the case $ d_1=1$.\\
By \eqref{Lem7-2b} and \eqref{Lem7-2}, we have
\begin{equation}  \nonumber 
        \fd_j: = \fa_{j,1} + \fa_{j,2} + \fa_{j,3} -m = \sum_{\ell=2}^5 \fa_{j,\ell} -m \in (-\mt_m, V_0), \quad j \in[1,4].
\end{equation}
From \eqref{Lem7-1} and \eqref{Lem7-2}, we get for $k_2 -k_1 > 10sV_0$ that $a_{1,1} =\fa_{1,1} \leq k_1$, $  \ff^{'}_{A_3,A_4,k_2}=0$. Hence $a_{3,1} =\fa_{3,1} =k_2$ and
\begin{equation}  \nonumber 
   10s V_0 \leq k_2 -k_1 = \fa_{3,1} -\fa_{1,1} = \fd_3 -\fd_1  \leq 2V_0, \quad V_0=10 \log_2 m.
\end{equation}
 We have a contradiction. \\

  Let us consider the cases $d_5=0$ and $d_2+d_3 =0$. If $d_5=0$, then $ 0= a_{3,5} - a_{1,5} \geq =k_2 -k_1 \geq 10sV_0$ (see \eqref{Lem7-2}). We have a contradiction. By the condition of Lemma 6, $d_1 <m$.
From \eqref{Lem4-5}, we get that if $d_2+d_3 =0$, then $m>d_1=d_1+d_2+d_3 =m$. We have a contradiction. \\

 Let us consider the case $d_1 \geq 2, \; d_5 \geq 1$ and $d_2+d_3 \geq 1$ :\\
By   Lemma C and  Lemma 5, we get
\begin{equation}   \nonumber 
  \sum_{j=1}^4 \fm^{[j]}=0. \qquad \fm^{[j]} \in   G_{A_j}.
\end{equation}

Let $\fm^{[5]} :=\fm^{[1]} + \fm^{[2]}$ and $\fm^{[6]}: =\fm^{[1]} + \fm^{[3]}$. We see that
 $\fm^{[5]} =-\fm^{[3]} - \fm^{[4]}$ and  $\fm^{[6]} =-\fm^{[2]} - \fm^{[4]}$.
 Hence $\rho(\fm_i^{[5]}) \leq \min(\rho(\fm_i^{[1]}+\fm_i^{[2]}),  \rho(\fm_i^{[3]}+\fm_i^{[4]})) $.
  In view of Lemma C, Lemma 5, \eqref{Lem7-2} and \eqref{Lem7-2a}, we obtain
\begin{equation} \label{Lem7-7}
     \rho(\fm_i^{[5]})  \leq
       \begin{cases}
   a_{2,i},  & \; \for  \; \; i \in J_2 \\
   a_{1,i},  & \; \for  \; \; i \in J_3 \\
   a_{3,i},  & \; \for  \; \; i \in \cJ_{4,\bA} \\
   a_{1,i},  & \; \for  \; \; i \in \cJ_{5,\bA}
  \end{cases} .
\end{equation}
Suppose that $\fm^{[5]}=0$  ($\fm{[1]}=\fm^{[2]}$, $\fm{[3]}=\fm^{[4]}$). Using  \eqref{Lem7-2} and \eqref{Lem7-2ab}, we get
\begin{multline} \label{Lem7-7ab}
a_{1,i} = a_{2,i}, \;\;\;a_{3,i} = a_{4,i} \;\;\;\;\; \for \;\; i \in J_1, \quad \quad
a_{1,i} = a_{3,i}, \;\;\;a_{2,i} = a_{4,i} \;\; \;\; \for \;\; i \in J_2 \cup J_3,\\
a_{1,i} - a_{2,i} \in [-4V_0,4V_0] \quad i=1,...,s. \quad
 \max_{j=1,2, i \in [1,s]} a_{j,i}=k_1, \quad  \max_{j=3,4, i \in [1,s]} a_{j,i}=k_2, \quad \\
   a_{j,1}+ \cdots +  a_{j,s}   \in ( m-\mt_m,m+V_0) \quad
\for \quad A_j \in I_{m,s,k}.
\end{multline}
Hence
\begin{equation} \label{Lem7-7abc}
  \Upsilon_{d_1,d_2,d_3,J,k} \ll    \sum_{\substack{ A_1,A_2 \in I_{m,s,k_1} \\ A_3,A_4 \in I_{m,s,k_2}}} \cL_{A,k} \d1(\eqref{Lem7-7ab} {\rm\; is\; true} ) \ll m^{2(s-2)-d_2-d_3}V_0^4 \ll
     m^{2(s-2)-9/10} .
\end{equation}

Now let $\fm_i^{[5]} \neq 0$.
  Bearing in mind that $ \fm^{[5]} \in \PP_m^{\bot}$, we get from \eqref{Lem7-7}
\begin{equation} \label{Lem7-8}
 \fa_{2,2} +  \fa_{1,3} +  \fa_{3,4} +  \fa_{1,5} \geq  \brho(\fm^{[5]}) \geq m -\mt_m  .
\end{equation}
Taking into account  that $\fa_{1,2}   + \fa_{1,3} +  \fa_{1,4} + \fa_{1,5} \leq m+V_0$ (see \eqref{Lem7-2b}), we obtain from \eqref{Lem7-8},  \eqref{Lem7-32} and \eqref{Lem7-36} that
\begin{multline} \label{Lem7-35}
 (\fa_{2,2}- \fa_{1,2}) +\fa_{1,2}  +  \fa_{1,3} +  \fa_{1,4} +  \fa_{1,5} -(\fa_{1,4} -\fa_{3,4}) =
 \fc_{1}- \fc_{2} +\fa_{1,2}  +  \fa_{1,3} +  \fa_{1,4} +  \fa_{1,5}
  \geq   m -\mt_m  ,\\
\fc_{1}- \fc_{2} \geq m -\mt_m -(\fa_{1,2}   + \fa_{1,3} +  \fa_{1,4} + \fa_{1,5}) \geq
  m -\mt_m -m -V_0 = -\mt_m-V_0.
\end{multline}\\

Let us consider $\fm^{[6]}= \fm^{[1]} + \fm^{[3]}$.
Suppose that $\fm_i^{[6]}=0$. Hence $\rho(\fm_i^{[1]}) =  \rho(\fm_i^{[3]}) $, $i=1,...,s$.
Applying \eqref{Lem7-2} and \eqref{Lem7-2ab}, we get
\begin{multline}  \nonumber 
 10s V_0 \leq k_2 -k_1 \leq \fa_{3,5} - \fa_{1,5} = \sum_{i \in \cJ_5} \big( \rho(\fm_i^{[3]}) -\rho(\fm_i^{[1]}) + (\epsilon_{3,i} - \epsilon_{1,i}) (V_0 +\mt_m) \big) \leq 4sV_0.
\end{multline}
We have a contradiction. So $\fm^{[6]}\neq 0$.

Bearing in mind
that  $\fm_i^{[6]} =-\fm_i^{[2]} - \fm_i^{[4]}$, we get  $\rho(\fm_i^{[6]}) \leq \min(\rho(\fm_i^{[1]}+\fm_i^{[3]}),  \rho(\fm_i^{[2]}+\fm_i^{[4]})) $.
Similarly, to \eqref{Lem7-7}, we have from \eqref{Lem7-2}  that
\begin{equation}  \nonumber 
     \rho(\fm_i^{[6]})  \leq
       \begin{cases}
   a_{1,i},  & \; \for  \; \; i \in J_2 \\
   a_{2,i},  & \; \for  \; \; i \in J_3 \\
   a_{1,i},  & \; \for  \; \; i \in \cJ_{4,\bA} \\
   a_{3,i},  & \; \for  \; \; i \in \cJ_{5,\bA}
  \end{cases} .
\end{equation}
Taking into account that $ \fm^{[6]} \neq 0$ and $ \fm^{[6]} \in \PP_m$, we get
\begin{equation} \label{Lem7-12}
  \fa_{1,2} +  \fa_{2,3} +  \fa_{1,4} +  \fa_{3,5} \geq  \brho(\fm^{[6]}) \geq m -\mt_m.
\end{equation}
 Bearing in mind  that $\fa_{1,2}   + \fa_{1,3} +  \fa_{1,4} + \fa_{1,5} \leq m+V_0$ (see \eqref{Lem7-2b}), we obtain from \eqref{Lem7-12} and    \eqref{Lem7-30} - \eqref{Lem7-36} that
\begin{multline} \nonumber
 (\fa_{2,3}- \fa_{1,3}) +\fa_{1,2}  +  \fa_{1,3} +  \fa_{1,4} +  \fa_{1,5} +(\fa_{3,5} -\fa_{1,5}) =
 \hat{\fc}_{2}- \hat{\fc}_{1} +\fa_{1,2}  +  \fa_{1,3} +  \fa_{1,4} +  \fa_{1,5}
  \geq   m -\mt_m  ,\\
\hat{\fc}_{2}- \hat{\fc}_{1} \geq m -\mt_m -(\fa_{1,2}   + \fa_{1,3} +  \fa_{1,4} + \fa_{1,5}) \geq
  m -\mt_m -m -V_0 = -\mt_m-V_0.
\end{multline}
Using \eqref{Lem7-30} - \eqref{Lem7-36}, we have $|\hat{\fc}_{j} - \fc_{j}| \leq 2V_0, \; j=1,2$, and
\begin{equation} \nonumber
   \fc_{2}- \fc_{1} \geq -6 V_0.
\end{equation}
By \eqref{Lem7-35}, we get
\begin{equation} \label{Lem7-70}
 \fa_{3,5} - \fa_{1,5} - \fa_{2,2} +  \fa_{1,2} =\fc_{2}- \fc_{1}  \in (-6V_0,2 V_0).
\end{equation}
From \eqref{Lem7-80}, we get
\begin{multline} \label{Lem7-87}
   \fg_{5,k_1,k_2} = \hat{\fg}_{5,k_1,k_2} + \check{\fg}_{5,k_1,k_2},  \quad \with \quad
 \hat{\fg}_{5,k_1,k_2} := |\Upsilon_{d_1,d_2,d_3,J,k}|  \d1( k_2 -k_1 > 10sV_0) \\
 \times  \d1(\fm^{[1]} =\fm^{[2]} ), \quad
 \check{\fg}_{5,k_1,k_2} := |\Upsilon_{d_1,d_2,d_3,J,k}|  \d1( k_2 -k_1 > 10sV_0)\d1(\fm^{[1]} \neq \fm^{[2]} ) .
\end{multline}
Applying Lemma 5, \eqref{Le3-5}, \eqref{Lem7-2} and \eqref{Lem7-70},   we obtain
\begin{multline} \nonumber
 \check{\fg}_{5,k_1,k_2} \ll    m^{2/5} \sum_{\substack{ A_1,A_2 \in I_{m,s,k_1} \\ A_3,A_4 \in I_{m,s,k_2}}}
\d1\big(  \fa_{3,5} - \fa_{1,5} - \fa_{2,2} +  \fa_{1,2}   \in (-6V_0,6V_0) \big)\;
W_1  W_2\cdots  W_5 \\
 \ll  m^{2/5} \sum_{A_j \in I_m^s, j =1,...,4} \; \prod_{j=1}^4 \Big(\d1( \max_{i} a_{j,i} =k_{[(j+1)/2]} ) \; \d1\big( \sum_{i=1}^s a_{j,i} \in (m-\mt_m, m+V_0)\big) \Big) \\
  \times \d1\big(  \fa_{3,5} - \fa_{1,5} - \fa_{2,2} +  \fa_{1,2}   \in (-V_0,V_0) \big) \\
 \times  \prod_{i\in \cJ_1} \d1(a_{1,i}=a_{2,i} ) \d1(a_{3,i}=a_{4,i} )
 \prod_{i\in \cJ_2 \cup \cJ_3} \d1(a_{1,i}=a_{3,i} ) \d1(a_{2,i}=a_{4,i} ).
\end{multline}
We fix $k_1,k_2$. We can choose $A_1$ in $O(m^{s-2} \log_2 m)$ ways. For given $A_1$,
we can choose $A_2$ in $O(m^{d_2+d_3-1} \log_2 m)$ ways.
For given $A_1,A_2$,
we can choose $A_3$ in $O(m^{d_1-1} \log_2 m)$ ways without  taking into account the ratio \eqref{Lem7-70}, and taking into account this ratio in $O(m^{d_1-2} \log_2 m)$  ways.
For given $A_1,A_2, A_3$, $A_4$ is chosen in the only one way. We have $d_1+d_2+d_3=s$.
Hence, we can choose $A_1,A_2, A_3, A_4$  in  $O(m^{2(s-2)-1} \log_2^2 m)$  ways.
Therefore
\begin{equation} \nonumber
   \check{\fg}_{5,k_1,k_2} = O(m^{2(s-2) -3/5} \log_2^2 m).
\end{equation}
Let us consider $ \check{\fg}_{5,k_1,k_2}  $. From \eqref{Lem7-7abc}, we get
\begin{equation} \nonumber
   \hat{\fg}_{5,k_1,k_2} = O(m^{2(s-2) -9/10}).
\end{equation}
By \eqref{Lem7-80} and \eqref{Lem7-87}, we have
\begin{equation} \nonumber
     \fg_{5} = O(m^{2(s-1) -2/5}).
\end{equation}
In view of \eqref{Lem7-85}, we get the assertion of the Lemma 6. \qed \\

By Lemma 4 - Lemma 6,  \eqref{Lem3-10a}  is proved. From \eqref{Lem3-8}, \eqref{Lem3-10} and  \eqref{Lem3-10a}, we get the bound for {\it Levi's conditional expectation} : \\ \\
{\bf Corollary 2.}
\begin{equation}  \nonumber 
 E_{Y,T}^{(m)} \Big(\Big(  \sum_{1 \leq k \leq m}
   \big( E_{Y,T}^{(m)} \big(\DD^2_k \; | \; \cF_{k-1} \big)
    -  E_{Y,T}^{(m)}( \DD_k^2) \big)   \Big)^2 \Big) \ll  m^{2(s-1) -1/5}.
\end{equation} \\ \\
{\bf 5. End of the proof of the theorems.}

Let
\begin{equation} \label{Lem8-1}
\dot{\SS}_m =\sum_{k=1}^m \DD_k  , \quad
     \dot{\varrho}_m^2 =E_{Y,T}^{(m)}(\dot{\SS}^2_m) \;\; \ad \;\;
\ddot{\SS}_m =D(\cP_m \oplus T,Y), \quad
     \ddot{\varrho}_m^2 =\EE_{Y,T,m} (\ddot{\SS}^2_m)    .
\end{equation}
By Lemma 1 and Lemma 2, we get
\begin{equation} \label{Lem8-1a}
	  \EE_{Y,T,m} ((\ddot{\SS}_m- \dot{\SS}_m)^2) =   \EE_{Y,T,m}((\cE^{(m)}(\cP_{m}\oplus T , Y ) +R(T,Y))^2) \ll  m^{1/5}   .
\end{equation}
From Lemma 3, we obtain
\begin{equation} \label{Lem8-1b}
\sum_{k=1}^m   E_{Y,T}^{(m)}(\DD_k^4)   \ll m^{2(s-1)-1/5} .
\end{equation}
In view of \eqref{In8}, we have for $\mt_m \leq 1/10 \;  \log_2 m $
\begin{equation} \label{Lem8-1c}
 2^{-4s-4} s^{-s+1}  \leq  	\ddot{\varrho}^2_m  m^{-s+1} = \cM_{s,2}^2 (\cP_{m})  m^{-s+1}   \leq  s^2 2^{2\mt_m+6}  \ll m^{1/5} .
\end{equation}
From (\ref{Lem8-1}), (\ref{Lem8-1a}),  (\ref{Lem8-1c}), (\ref{Le2-5}) and (\ref{Le2-5a}), we obtain
\begin{multline}   \label{Lem8-3}
   \dot{\varrho}_m^2= E_{Y,T}^{(m)}(\dot{\SS}^2_m) = \EE_{Y,T,m}(\dot{\SS}^2_m)
   =\EE_{Y,T,m}((\ddot{\SS}_m -(\dot{\SS}_m-\dot{\SS}_m))^2  ) \\
    \geq 1/2 \EE_{Y,T,m}((\ddot{\SS}_m)^2) -  \EE_{Y,T,m}\big((\dot{\SS}_m-\dot{\SS}_m)^2\big)
  \geq 1/2 \ddot{\varrho}_m^2 \\
	- O(m^{2(s-1)-1/5} )     \geq 2^{-2s-4} s^{-s+1}  m^{s-1}
\end{multline}
for $m \geq m_0$ with some $m_0 >0$.

By \eqref{Lem8-1a} and \eqref{Lem8-1c}, we get
\begin{multline}   \label{Lem8-4}
|\ddot{\varrho}_m^2 -\dot{\varrho}_m^2|^2 = \big|\EE_{Y,T,m}^2 \big(  \ddot{\SS}_m^2 - \dot{\SS}_m^2\big) \big| =   \EE_{Y,T,m}^2\big|\big(  \ddot{\SS}_m - \dot{\SS}_m \big)
   \big(  \ddot{\SS}_m + \dot{\SS}_m \big) \big| \leq
   \EE_{Y,T,m} \big(  (\ddot{\SS}_m - \dot{\SS}_m )^2 \big)  \\
  \times   \EE_{Y,T,m}  ( (\ddot{\SS}_m +\dot{\SS}_m )^2 )
  \ll  m^{1/5}
   \EE_{Y,T,m}  \Big(\big( 2\ddot{\SS}_m -(\ddot{\SS}_m -\dot{\SS}_m)      \big)^2 \Big) \\
\ll    m^{1/5}  \Big( \EE_{Y,T,m} \big( \ddot{\SS}_m^2 \big) +
  \EE_{Y,T,m} \big((\ddot{\SS}_m -\dot{\SS}_m)^2      \big) \Big)   \ll m^{s-1 +2/5} .
\end{multline}\\

 We shall use  the following variant of the {\it martingale central limit theorem} (see [Ha, p.~58,
 Corollary 3.1]):

 Let $(\Omega, \cF,P)$ be a probability space and $\{(\zeta_{m,k}, \FF_{m,k}) \;| \;  k=1,...,\ell_m\}$ be a martingale difference array with
 ${\bf E} (\zeta_{m,k} \;|\; \FF_{m,k-1}) =0$ a.s. ($\FF_{m,0}$ is the trivial field).
 \\ \\
{\bf Theorem  A.} {\it Let
\begin{multline}\label{Le13-10}
 \SS_{m} =\sum_{1 \leq k \leq \ell_m} \zeta_{m,k}, \quad
 L(m,\epsilon) = \sum_{1 \leq k \leq \ell_m} {\bf E} (\zeta_{m,k}^2 \d1(|\zeta_{m,k}|>\epsilon)),
  \quad
  \sum_{1 \leq k \leq \ell_m}  {\bf E} (\zeta_{m,k}^2) =1,\\
  \VV_{m}^2 =\sum_{1 \leq k \leq \ell_m} {\bf E} (\zeta_{m,k}^2 \; | \; \FF_{m,k-1}) ,  \quad
  L(m,\epsilon)  \stackrel{P}{\rightarrow} 0
    \quad \forall  \epsilon>0,  \quad   \VV_{m}^2 \stackrel{P}{\rightarrow} 1.
\end{multline}
Then $ \SS_{m} 	\stackrel{w}{\rightarrow} \cN(0,1)    $. } \\

 By (\ref{Lem3-7}),  $\{(\DD_k, F_{k}) \;| \;
     k=1,...,m\}$ is the martingale difference array.

   Now we apply Theorem A to  the  array $\{(\DD_k, F_{k}) \;| \;
     k=1,...,m\}$ with $\FF_{m,k}= F_k$,
 $\zeta_{m,k} = \DD_k/\dot{\varrho}_m$
 and  $\ell_m= m$.\\ \\
{\bf Lemma 7.} {\it With the notations as above},
   $ \dot{\SS}_m / \dot{\varrho}_m 	\stackrel{w}{\rightarrow} \cN(0,1)  $. \\ \\
{\bf Proof.}
Let us consider $ L(m,\epsilon)$. Using   \eqref{Lem8-1b}, \eqref{Lem8-3} and \eqref{Le13-10},  we get
\begin{multline}\nonumber
  L(m,\epsilon) = \sum_{1 \leq k \leq m}
   E_{Y,T}^{(m)}  \big(  |\DD_k / \dot{\varrho}_m|^{2} \d1 (|\DD_k/\dot{\varrho}_m|> \epsilon)  \Big)  \\
   \leq \sum_{1 \leq k \leq m}    E_{Y,T}^{(m)}  \Big(   |\DD_k / \dot{\varrho}_m|^{2}\;  \frac{|\DD_k/\dot{\varrho}_m|^2}{\epsilon^{2} } \Big)
   =  \sum_{1 \leq k \leq m}\epsilon^{-2}
   E_{Y,T}^{(m)}  \big((  \DD_k / \dot{\varrho}_m)^{4}  \big)\\
 \ll  \epsilon^{-2} \sum_{1 \leq k \leq m}  m^{-2(s-1)} \;
   \big(E_{Y,T}^{(m)}  (\DD^4_k)) \big)   \ll \epsilon^{-2}   m^{-2(s-1)}  m^{2(s-1)-1/5}
     \ll m^{-1/5} .
\end{multline}
Let us consider $\VV_m^2 -1$.
By \eqref{Le13-10}, Corollary 2  and  Chebyshev's inequality, we have
\begin{multline*} 
  P( |\VV_m^2 -1|  >\epsilon)   \leq  \epsilon^{-2}  \Big( E_{Y,T}^{(m)} \big( |\VV_m^2 -1|^2 \big) \\
 = E_{Y,T}^{(m)} \Big(\Big(  \sum_{1 \leq k \leq m}  \dot{\varrho}_m^{-2}
   \big( E_{Y,T}^{(m)} \big(\DD^2_k \; | \; F_{k-1} \big)
    -  E_{Y,T}^{(m)} (\DD_k^2) \big)   \Big)^2 \Big)\\
           \ll   m^{-2(s-1)}  m^{2(s-1)-1/5}      \ll m^{-1/5} .
\end{multline*}
Applying Theorem A, we obtain the assertion of Lemma 7. \qed \\ \\
%
We need the following ``Converging Together Lemma'' : \\ \\
{\bf Lemma E.} [Du, p.105, ex.3.2.13] {\it If $U_m \stackrel{w}{\rightarrow} U$ and
$Z_m - U_m \stackrel{w}{\rightarrow} 0$,
then $ Z_m \stackrel{w}{\rightarrow} U$.}\\ \\
{\bf  Proof of Theorem 1.}
In view of  \eqref{Lem8-1} - \eqref{Lem8-4}, we get
\begin{multline}  \nonumber 
\dot{\varrho}_m^2 \geq 2^{-4s-10} s^{-s}m^{s-1}  , \quad \ddot{\varrho}_m \geq  2^{-4s-8} s^{-s} m^{s-1},  \quad |\dot{\varrho}_m^2- \ddot{\varrho}_m^2|^2 \ll m^{s-3/5}  ,\\
 \EE_{Y,T,m} ( \dot{\SS}_m -\ddot{\SS}_m )^2 \ll m^{1/5}  \quad \ad \quad
 \EE_{Y,T,m} ( \ddot{\SS}_m )^2 \ll m^{s-1+1/5}.
\end{multline}
and
\begin{equation} \nonumber 
    \Big| \frac{1}{\dot{\varrho}_m} -   \frac{1}{\ddot{\varrho}_m} \Big|^2
      =
        \frac{ | \dot{\varrho}_m - \ddot{\varrho}_m |^2}{\dot{\varrho}_m^2  \ddot{\varrho}_m^2}
     =
       \frac{ |\dot{\varrho}_m^2- \ddot{\varrho}_m^2|^2 } {\dot{\varrho}_m^2
       \ddot{\varrho}_m^2 (\dot{\varrho}_m+ \ddot{\varrho}_m)^2} \ll m^{s-3/5-3(s-1)}=
       m^{-2s+12/5} .
\end{equation}
%
By  \eqref{Lem8-1} - \eqref{Lem8-4}, we obtain
\begin{multline*}
 \EE_{Y,T,m} \Big( \frac{ \dot{\SS}_m}{\dot{\varrho}_m}  - \frac{ \ddot{\SS}_m}{\ddot{\varrho}_m} \Big)^2 =
 \EE_{Y,T,m} \Big( \frac{ \dot{\SS}_m -\ddot{\SS}_m }{\dot{\varrho}_m} +
  \ddot{\SS}_m ( \frac{1}{\dot{\varrho}_m}  -\frac{1}{\ddot{\varrho}_m})\Big)^2
     \leq 2 \dot{\varrho}^{-2}_m   \EE_{Y,T,m} ( \dot{\SS}_m -\ddot{\SS}_m )^2  \\
     +2  (1/\dot{\varrho}_m  -1/\ddot{\varrho}_m)^2  \EE_{Y,T,m}  \ddot{\SS}^2_m  \\
     \ll    m^{1/5 -(s-1)} +
     m^{s-1 +1/5 +(-2s+12/5)}  \ll  m^{-s +8/5} \;\stackrel{m \to \infty}{\longrightarrow} 0  \quad \for \; s \geq 2 .
\end{multline*}
Hence
$\dot{\SS}_m /\dot{\varrho}_m - \ddot{\SS}_m /\dot{\varrho}_m	\stackrel{w}{\rightarrow} 0$. Bearing in mind that
$ \dot{\SS}_m /\dot{\varrho}_m	\stackrel{w}{\rightarrow} \cN(0,1)  $ and
  Lemma~E, we get that
$\ddot{\SS}_m /\ddot{\varrho}_m	\stackrel{w}{\rightarrow} \cN(0,1)  $.
 By \eqref{Lem8-1}, we have  the assertion of Theorem 1. \qed \\ \\

{\bf  Proof of Theorem 2.}  We need the following simple variant of the
 {\it Continuous Mapping Theorem} (see [Du, Theorem 3.2.4., p.101]).\\ \\
{\bf Theorem B.} {\it Let $g$ be a continuous function. If $U_m \stackrel{w}{\rightarrow}   U$, then } $g(U_m) \stackrel{w}{\rightarrow}   g(U)$.

By [Bil, p.31 ], a simple condition of uniform integrability of a sequence of functions
 $U_m$ is that $\sup_m E|U_m|^{1+\epsilon} < \infty$. According to [Bil, Theorem 3.5, p.31], we have\\ \\
{\bf Theorem C.}
 {\it If $U_m$ are uniformly integrable and $U_m \stackrel{w}{\rightarrow}   U$,
  then $U$ is integrable and $E(U_m) \to E(U)$ .   } \\ \\
Let  $ Z_m:=  D( \cP_m \oplus T,Y) / \cM_{s,2} (\cP_{m}) $.
By Theorem 2,  $Z_m \stackrel{w}{\rightarrow}   \cN(0,1)=:Z$. We take the continuous function $g(x)=|x|^p$.
Using Theorem B, we get $g(Z_m) \stackrel{w}{\rightarrow}   g(Z)$.
   Bearing in mind \eqref{In8}, we get that the functions $g(Z_m)\; (m=1,2,...)$ are uniformly integrable.

   Now using Theorem C, we get the assertion of Theorem 2. \qed \\ \\

{\bf Address}: Department of Mathematics, Bar-Ilan University, Ramat-Gan, 5290002, Israel \\
{\bf E-mail}: mlevin@math.biu.ac.il\\
\end{document}